# The needlets bispectrum[*]

## Xiaohong Lan

*Institute of Mathematics, Chinese Academy of Sciences and
Department of Mathematics, University of Rome Tor Vergata
e-mail:* lan@mat.uniroma2.it

## Domenico Marinucci

*Department of Mathematics, University of Rome Tor Vergata
e-mail:* marinucc@mat.uniroma2.it

**Abstract:** The purpose of this paper is to join two different threads of the recent literature on random fields on the sphere, namely the statistical analysis of higher order angular power spectra on one hand, and the construction of second-generation wavelets on the sphere on the other. To this aim, we introduce the needlets bispectrum and we derive a number of convergence results. Here, the limit theory is developed in the high resolution sense. The leading motivation of these results is the need for statistical procedures for searching non-Gaussianity in Cosmic Microwave Background radiation.

Received February 2008.

## Contents

[*]We are grateful to a Referee and an Associate Editor for a very careful reading and many insightful remarks on an earlier version of this paper. We are also grateful to P.Baldi, G.Kerkyacharian and D.Picard for many useful discussions. The research of Xiaohong Lan was supported by the Equal Opportunity Committee of the University of Rome Tor Vergata.

## 1. Introduction

The purpose of this paper is to join two different threads of the recent literature on random fields on the sphere, namely the statistical analysis of higher order angular power spectra on one hand, and the construction of second-generation wavelets on the sphere on the other hand. More precisely, we shall be concerned with zero-mean, mean square continuous and isotropic random fields on the sphere, for which the following spectral representation holds [1, 4]: for $x \in S^2$,

$$T(x) = \sum_{l,m} a_{lm} Y_{lm}(x) \tag{1}$$

where $\{a_{lm}\}_{l,m}$, $m = 1, \ldots, l$ a triangular array of zero-mean, orthogonal, complex-valued random variables with variance $\mathbb{E}|a_{lm}|^2 = C_l$, the angular power spectrum of the random field. For $m < 0$ we have $a_{lm} = (-1)^m \overline{a_{l-m}}$, whereas $a_{l0}$ is real with the same mean and variance. (1) holds in the $L^2(S^2)$ sense, i.e. we have

$$\lim_{L \to \infty} E \left| \int_{S^2} \left\{ T(x) - \sum_{l=1}^{L} \sum_{m=-l}^{l} a_{lm} Y_{lm}(x) \right\}^2 \mu(dx) \right|^2 = 0,$$

$\mu(dx)$ denoting the uniform measure on the sphere. The functions $\{Y_{lm}(x)\}$ are the so-called spherical harmonics, i.e. the eigenvectors of the Laplacian operator on the sphere,

$$
\begin{aligned}
\Delta_{S^2} Y_{lm}(\vartheta, \varphi) &= \left[ \frac{1}{\sin \vartheta} \frac{\partial}{\partial \vartheta} \left\{ \sin \vartheta \frac{\partial}{\partial \vartheta} \right\} + \frac{1}{\sin^2 \vartheta} \frac{\partial^2}{\partial \varphi^2} \right] Y_{lm}(\vartheta, \varphi) \\
&= -l(l+1) Y_{lm}(\vartheta, \varphi)
\end{aligned}
$$

where we have moved to spherical coordinates $x = (\vartheta, \varphi)$, $0 \leq \vartheta \leq \pi$ and $0 \leq \varphi < 2\pi$. It is a well-known result that the spherical harmonics provide a complete orthonormal systems for $L^2(S^2)$ [36].

The analysis of random fields on the sphere has recently gained very strong physical motivations, due to the overwhelming amount of data which is becoming available on Cosmic Microwave Background radiation (hereafter (CMB)). As detailed elsewhere [15], to the first order we can view CMB data as maps of the Universe in the immediate adjacency of the Big Bang. The first of these maps were provided by the satellite experiment *COBE* in 1993, and in view of this G. Smoot and J. Mather were awarded the Nobel prize for Physics in 2006. Much more refined maps have been made available by another NASA satellite experiments, *WMAP* (http://wmap.gsfc.nasa.gov/); still more refined data are expected from the ESA mission *Planck,* scheduled to be launched in October 2008. These huge collaborations involve hundred of scientists and are expected to provide invaluable information on Physics and Cosmology. At the same time, these massive data sets have called for huge statistical challenges, ranging from power spectrum estimation to outlier detection, testing for isotropy, efficient denoising and map-making, handling with missing data, testing for non-Gaussianity and many others.

Among these issues, particular interest has been driven by efficient testing for non-Gaussianity. This is due to strong physical motivations (the leading paradigm for the Big Bang dynamics predicts (very close to) Gaussian fluctuations) and difficulties in finding a proper statistical procedure. There is by now a wide consensus that the most efficient procedures to probe non-Gaussianity are based upon the bispectrum, in the idealistic circumstances where the spherical random field is fully observed: see for instance [8, 11, 20, 23, 24, 37]. However, the properties of the bispectrum are also known to deteriorate dramatically in the presence of missing observations, see again [11]. To handle the latter problem, a general approach is to focus on wavelet, rather than Fourier transforms. The construction of spherical wavelets has recently drawn an enormous amount of attention in the literature, see for instance [2, 27, 38] and the references therein. In our view, a particularly convenient tight frame construction on the sphere is provided by so-called needlets, which were introduced in [30, 31]; their applications to spherical random fields is due to [6, 7]. Needlets enjoy two properties which seem especially worth recalling: they are quasi-exponentially localized in the real domain and compactly supported in the harmonic domain. A further, quite unexpected property is as follows: random needlets coefficients are asymptotically uncorrelated at the highest frequencies and hence, in the Gaussian case, independent, see again [6]. This latter feature is rather surprising in a compact domain and makes asymptotic theory possible even in the presence of a single realization of a spherical random field.

Our aim in this paper is to borrow ideas from the bispectrum and the needlets literature to propose and analyze a needlets bispectrum, where the random coefficients in the needlets expansion are combined in a similar way to the bispectrum construction. The aim is to obtain a procedure which mimicks the ability of the bispectrum to search for non-Gaussianity at the most efficient combination of frequencies, at the same time providing a much more robust construction in the presence of missing data, as typical of the needlets. The plan of the paper is as follows: in Section 2 we review some background material on spherical random fields, the bispectrum and the needlets construction. In Section 3 we introduce the needlets bispectrum and we establish a central limit theorem, in the high resolution sense. In Section 4 we go on to establish a functional central limit theorem for the integrated needlets bispectrum; in Section 5 we provide some preliminary discussion on the behaviour under non-Gaussian assumptions and discuss some possibilities for applications and further research.

## 2. Some background material

### 2.1. A Review on Needlets

In this subsection we shall briefly recall the main features of the needlets construction. As mentioned above, needlets were first introduced in the Functional

Analysis literature by [30, 31], whereas the investigation of their properties from a stochastic point of view is due to [5, 6] and [7]; see also [19, 17]. We need first to introduce some notation and definitions, which are largely identical to those in [6, 31].

Given any two positive sequences $\{a_j\}$, $\{b_j\}$, we write $a_j \approx b_j$ if there exist $c > 0$ such that $c^{-1}a_j \leq b_j \leq ca_j$ for all $j$. The standard (open and closed) balls in $\mathbb{S}^2$ are given as always by $B(a, \alpha) = \{x, \ d(a, x) \leq \alpha\}$, $B^\circ(a, \alpha) = \{x, \ d(a, x) < \alpha\}$. For a general $A \subset \mathbb{S}^2$ we will denote by $|A|$ the spherical measure of $A$. Now fix $\epsilon > 0$ and $x_1, \ldots x_N$ in $\mathbb{S}^2$ such that $\forall i \neq j, \ d(x_i, x_j) > \epsilon$; the set $\{x_1, \ldots, x_N\} = \Xi_\epsilon$ is called a maximal $\epsilon-$ net if it satisfies $\forall x \in \mathbb{S}^2, \ d(x, \Xi_\epsilon) \leq \epsilon, \ \cup_{x_i \in \Xi_\epsilon} B(x_i, \epsilon) = \mathbb{S}^2$ and $\forall i \neq j, \ \ B(x_i, \epsilon/2) \cap B(x_j, \epsilon/2) = \emptyset$. It follows from Lemma 5 in [7] that

$$\frac{4}{\epsilon^2} \leq N \leq \frac{4}{\epsilon^2}\pi^2 \tag{2}$$

For all $x_i \in \Xi_\epsilon$, the associated family of *Voronoi cells* is defined by:

$$\mathcal{V}(x_i) = \{x \in \mathbb{S}^2, \ \ \forall j \neq i, \ \ d(x, x_i) \leq d(x, x_j)\}.$$

We recall that $B(x_i, \epsilon/2) \subset \mathcal{V}(x_i) \subset B(x_i, \epsilon)$. Also, if two Voronoi cells are adjacent, i.e. $\mathcal{V}(x_i) \cap \mathcal{V}(x_j) \neq \emptyset$, then by necessity $d(x_i, x_j) \leq 2\epsilon$. It is proved in [7] that there are at most $6\pi^2$ adjacent cells to any given cell.

For the construction of needlets, we should first start to define $\mathcal{K}_l$ as the space of the restrictions to the sphere $\mathbb{S}^2$ of polynomials of degree less than $l$. The next ingredient are the set of *cubature points* and *cubature weights;* indeed, it is now a standard result (see for instance [29]) that for all $j \in \mathbb{N}$, there exists a finite subset $\mathcal{X}_j$ of $\mathbb{S}^2$ and positive real numbers $\lambda_{jk} > 0$, indexed by the elements of $\mathcal{X}_j$, such that

$$\forall f \in \mathcal{K}_l, \ \ \int_{\mathbb{S}^2} f(x)dx = \sum_{\xi_{jk} \in \mathcal{X}_j} \lambda_{jk} f(\xi_{jk}). \tag{3}$$

Given a fixed $B > 1$, we shall denote by $\{\xi_{jk}\}$ the cubature points corresponding to the space $\mathcal{K}_{[3B^{j+1}]}$, where $[.]$ represents as usual integer part. It is known that $\{\mathcal{X}_j\}_{j=0}^\infty$ can be taken s.t. the cubature points for each $j$ are almost uniformly $\epsilon_j-$distributed with $\epsilon_j := \kappa B^{-j}$, and the coefficients $\{\lambda_{jk}\}$ are such that $\lambda_{jk} \approx cB^{-2j}$, $card\{\mathcal{X}_j\} \approx B^{2j}$.

Now let $\phi$ be a $C^\infty$ function supported in $|\xi| \leq 1$, such that $0 \leq \phi(\xi) \leq 1$ and $\phi(\xi) = 1$ if $|\xi| \leq B^{-1}$, $B > 1$. Following again [30, 31], we define

$$b^2(\xi) = \phi\left(\frac{\xi}{B}\right) - \phi(\xi) \geq 0 \text{ so that } \forall |\xi| \geq 1, \sum_j b^2\left(\frac{\xi}{B^j}\right) = 1. \tag{4}$$

It is immediate to verify that $b(\xi) \neq 0$ only if $\frac{1}{B} \leq |\xi| \leq B$. The needlets frame $\{\varphi_{jk}\}$ is then constructed as

$$\varphi_{jk}(x) := \sqrt{\lambda_{jk}} \sum_l b\left(\frac{l}{B^j}\right) \sum_{m=-l}^l Y_{lm}(\xi_{jk})Y_{lm}^*(x). \tag{5}$$

The main localization property of $\{\varphi_{jk}\}$ is established in [30], where it is shown that for any $M \in \mathbb{N}$ there exists a constant $c_M > 0$ s.t., for every $\xi \in \mathbb{S}^2$:

$$|\varphi_{jk}(\xi)| \leq \frac{c_M B^j}{(1 + B^j \arccos\langle \xi_{jk}, \xi \rangle)^M} \text{ uniformly in } (j, k).$$

More explicitly, needlets are almost exponentially localized around any cubature point, which motivates their name. In the stochastic case, the (random) spherical needlet coefficients are then defined as

$$\beta_{jk} = \int_{\mathbb{S}^2} T(x)\varphi_{jk}(x)dx = \sqrt{\lambda_{jk}} \sum_l b\left(\frac{l}{B^j}\right) \sum_{m=-l}^{l} a_{lm} Y_{lm}(\xi_{jk}). \quad (6)$$

We have immediately

$$\sum_k \beta_{jk} \sqrt{\lambda_{jk}} = 0, \quad (7)$$

i.e., the (weighted) sample mean of the needlets coefficients is identically zero at all levels $j$. The proof is trivial, because

$$
\begin{aligned}
\sum_k \beta_{jk} \sqrt{\lambda_{jk}} &= \sum_{l=B^{j-1}}^{B^{j+1}} \sum_{m=-l}^{l} b\left(\frac{l}{B^j}\right) a_{lm} \left[\sum_k \lambda_{jk} Y_{lm}(\xi_{jk})\right] \\
&= \sum_{l=B^{j-1}}^{B^{j+1}} \sum_{m=-l}^{l} b\left(\frac{l}{B^j}\right) a_{lm} \left[\int_{\mathbb{S}^2} Y_{lm}(x)dx\right] = 0.
\end{aligned}
$$

The variance of the needlets coefficients is given by

$$
\begin{aligned}
\mathbb{E}\beta_{jk}^2 &= \lambda_{jk} \sum_{l=B^{j-1}}^{B^{j+1}} b^2\left(\frac{l}{B^j}\right) C_l \sum_{m=-l}^{l} Y_{lm}(\xi_{jk})\overline{Y_{lm}(\xi_{jk})} \\
&= \lambda_{jk} \sum_{l=B^{j-1}}^{B^{j+1}} b^2\left(\frac{l}{B^j}\right) C_l \frac{2l+1}{4\pi} P_l(\cos 0) \\
&= \lambda_{jk} \sum_{l=B^{j-1}}^{B^{j+1}} \frac{2l+1}{4\pi} b^2\left(\frac{l}{B^j}\right) C_l =: \sigma_{jk}^2 > 0.
\end{aligned}
$$

Note that we have $\sigma_{jk}^2 \approx: \sigma_j^2$ uniformly over $k$, where

$$\sigma_j^2 := \frac{4\pi}{N_j} \sum_{l=B^{j-1}}^{B^{j+1}} \frac{2l+1}{4\pi} b^2\left(\frac{l}{B^j}\right) C_l, \ N_j = card\{\mathcal{X}_j\}.$$

From now on, we shall typically focus on the normalized needlets coefficients, defined as $\widehat{\beta}_{jk} := \beta_{jk}/\sigma_j$.

To investigate the correlation, we introduce now the same, mild regularity conditions on the angular power spectrum $C_l$ of the random field $T(x)$ as in [5, 6].

**Condition A** The random field $T(x)$ is Gaussian and isotropic with angular power spectrum such that, for all $B > 1$, there exist $\alpha > 2$, and $\{g_j(.)\}_{j=1,2,...}$ a sequence of functions such that

$$C_l = l^{-\alpha} g_j\left(\frac{l}{B^j}\right) > 0, \text{ for } B^{j-1} < l < B^{j+1}, \ j = 1, 2, \dots \tag{8}$$

where

$$c_0^{-1} \leq g_j \leq c_0 \text{ for all } j \in \mathbb{N}, \text{ and}$$

$$\sup_j \sup_{B^{-1} \leq u \leq B} \left| \frac{d^r}{du^r} g_j(u) \right| \leq c_r \text{ some } c_0, c_1, \dots c_M > 0, \ M \in \mathbb{N}.$$

Condition 8 entails a weak smoothness requirement on the behaviour of the angular power spectrum, which is trivially satisfied by some cosmologically relevant models (where the angular power spectrum usually behaves as an inverse polynomial, see again [15], pp. 243–244). For instance, considering positive constants $d_j, \ j = 1, \dots, p, \ p > 2$, for all $B > 1$ condition A holds for

$$
\begin{aligned}
C_l &= \frac{1}{\sum_{k=0}^p d_k l^k} = l^{-p} \frac{l^p}{d_0 + d_1 l + \dots d_p l^p} \\
&= l^{-p} \frac{(l/B^j)^p}{\sum_{k=0}^p d_k B^{jk-jp}(l/B^j)^k} \\
&= l^{-p} g_j\left(\frac{l}{B^j}\right), \text{ for } B^{j-1} < l < B^{j+1}, \ j = 1, 2, \dots \\
g_j(u) &:= \frac{u^p}{\sum_{k=0}^p d_k B^{jk-jp} u^k}.
\end{aligned}
$$

Under Condition 8 a crucial and rather unexpected property of the random needlets coefficients is established in [6], namely the correlation bound

$$|Corr(\beta_{jk}, \beta_{jk'})| \leq \frac{c_M}{(1 + B^j d(\xi_{jk}, \xi_{jk'}))^M} \tag{9}$$

where $d(\xi_{jk}, \xi_{jk'}) = \arccos(\langle \xi_{jk}, \xi_{jk'} \rangle)$. In words, (9) is stating that as the frequency increases, spherical needlets coefficients are asymptotically uncorrelated (and hence, in the Gaussian case, independent), for any given angular distance. This property is of course of the greatest importance when investigating the asymptotic behaviour of statistical procedures: in some sense, it states that it is possible to derive an infinitely growing array of asymptotically independent "observations" (the needlets coefficients) out of a single realization of a continuous random field on a compact domain. It should be stressed that this property is not by any means a consequence of the localization properties of the needlets frame. As a counterexample, it is easy to construct spherical frames having bounded support in real space, whereas the corresponding random coefficients are not at all uncorrelated (recall the angular correlation function can be taken to be bounded from below at any distance on the sphere).

We recall briefly two further results that we shall exploit often in the sequel; more precisely, ([6], *Lemma 10*) we note that for $M > 2$, $j \in \mathbb{N}$, $\forall$ $k, k'$,

$$\sum_{\xi_{jk} \in \chi_j} \frac{1}{(1 + B^j d(\xi_{jk}, \xi_{jk'}))^M} \leq C_M \text{ some } C_M > 0. \tag{10}$$

Also ([31], *Lemma 4.8*), for some $C'_M$ depending only on $M$ we have the inequality

$$\sum_{\xi_{jk} \in \chi_j} \frac{1}{(1 + B^j d(\xi_{jk}, \xi_{jk'}))^M} \frac{1}{(1 + B^j d(\xi_{jk}, \xi_{jk''}))^M} \leq \frac{C'_M}{(1 + B^j d(\xi_{jk'}, \xi_{jk''}))^M}. \tag{11}$$

From the computational point of view, we should stress that needlets are not only feasible, but indeed extremely convenient. The implementation can be performed with a minimal effort by means of standard packages for the analysis of spherical random fields such as HEALPIX or GLESP ([18] and [16]), see for details [26], where plots and numerical evidence on localization and uncorrelation are also provided.

## 2.2. Diagram Formula

To complete our background, we need a quick review on the diagram formula. This is a material which can now be found at a textbook level (see for instance [34]); nevertheless, we need a brief overview to fix notation. Denote by $H_q$ the $q-$th order Hermite polynomials, defined as

$$H_q(u) = (-1)^q e^{u^2/2} \frac{d^q}{du} e^{-u^2/2}.$$

We now introduce *diagrams*, which are basically mnemonic devices for computing the moments and cumulants of polynomial forms in Gaussian random variables. Our notation is the same as for instance in [23, 24]. Let $p$ and $l_j$, $j = 1, \ldots, p$ be given integers. A diagram $\gamma$ of order $(l_1, \ldots, l_p)$ is a set of points $\{(j, l) : 1 \leq j \leq p, 1 \leq l \leq l_j\}$ called vertices, viewed as a table $W = \overrightarrow{l_1} \otimes \cdots \otimes \overrightarrow{l_p}$ and a partition of these points into pairs

$$\{((j, l), (k, s)) : 1 \leq j \leq k \leq p; 1 \leq l \leq l_j, 1 \leq s \leq l_k\},$$

called edges. We denote by $\Gamma(W)$ the set of diagrams of order $(l_1, \ldots, l_p)$. If the order is $l_1 = \cdots = l_p = q$, for simplicity, we also write $\Gamma(p, q)$ instead of $\Gamma(W)$. We say that:

*a)* A diagram has a flat edge if there is at least one pair $\{(i, j)(i', j')\}$ such that $i = i'$; we write $\Gamma_F$ for the set of diagrams that has at least one flat edge, and $\Gamma_{\overline{F}}$ otherwise.

*b)* A diagram is *connected* if it is not possible to partition the rows $\overrightarrow{l_1} \cdots \overrightarrow{l_p}$ of the table $W$ into two parts, i.e. one cannot find a partition $K_1 \cup K_2 = \{1, \ldots, p\}$ that, for each member $V_k$ of the set of edges $(V_1, \ldots, V_r)$ in a diagram $\gamma$, either $V_k \in \cup_{j \in K_1} \overrightarrow{l_j}$, or $V_k \in \cup_{j \in K_2} \overrightarrow{l_j}$ holds; we write $\Gamma_C$ for connected diagrams, and $\Gamma_{\overline{C}}$ otherwise.

*c)* A diagram is paired if, considering any two sets of edges $\{(i_1, j_1)(i_2, j_2)\}$ $\{(i_3, j_3)(i_4, j_4)\}$, then $i_1 = i_3$ implies $i_2 = i_4$; in words, the rows are completely coupled two by two. We write $\Gamma_P$ for the set of diagrams for paired diagrams, and $\Gamma_{\overline{P}}$ otherwise.

**Proposition 1.** (Diagram Formula) *Let* $(z_1, \ldots, z_p)$ *be a centered Gaussian vector, and let* $\gamma_{ij} = E[z_i z_j], i, j = 1, \ldots, p$ *their covariances, Let* $H_{l_1}, \ldots, H_{l_p}$ *be Hermite polynomials of degree* $l_1, \ldots, l_p$ *respectively. Let* $L$ *be a table consisting of* $p$ *rows* $l_1, \ldots l_p$, *where* $l_j$ *is the order of Hermite polynomial in the variable* $z_j$. *Then*

$$E[\Pi_{j=1}^p H_{l_j}(z_j)] = \sum_{G \in \Gamma(l_1, \ldots, l_p)} \Pi_{1 \le i \le j \le p} \gamma_{ij}^{\eta_{ij}(G)}$$

$$Cum(H_{l_1}(z_1), \ldots, H_{l_p}(z_p)) = \sum_{G \in \Gamma_c(l_1, \ldots, l_p)} \Pi_{1 \le i \le j \le p} \gamma_{ij}^{\eta_{ij}(G)}$$

*where, for each diagram* $G$, $\eta_{ij}(G)$ *is the number of edges between rows* $l_i, l_j$ *and* $Cum(H_{l_1}(z_1), \ldots, H_{l_p}(z_p))$ *denotes the* $p-th$ *order cumulant.*

We have now all the preliminary material to define our needlets bispectrum on $\mathbb{S}^2$, as explained in the following Section.

## 3. A Central Limit Theorem for the Needlets Bispectrum

### 3.1. The needlets bispectrum

As mentioned in the introduction, the recent literature suggests that the most powerful statistic to search for non-Gaussianity in fully observed spherical random fields is the (normalized) angular bispectrum, defined as

$$I_{l_1 l_2 l_3} = \sum_{m_1 m_2 m_3} \begin{pmatrix} l_1 & l_2 & l_3 \\ m_1 & m_2 & m_3 \end{pmatrix} \frac{a_{l_1 m_1} a_{l_2 m_2} a_{l_3 m_3}}{\sqrt{C_{l_1} C_{l_2} C_{l_3}}}$$

where the symbol in brackets represents the so-called Wigner's $3j$ coefficients, which are meant to ensure the statistics is rotationally invariant. Wigner's $3j$ coefficients arise in many different instances, especially in the quantum theory of angular momentum (see [36], where explicit expressions are also provided). Up to a normalization factors, they are equivalent to the so-called Clebsch-Gordan coefficients, which play an important role in representation theory for the group of rotations $SO(3)$, see [25] for a much more detailed discussion and probabilistic applications. Following [20, 37], an alternative definition of the (normalized) bispectrum can be considered, namely

$$\widetilde{I}_{l_1 l_2 l_3} = \int_{S^2} \frac{T_{l_1}(x) T_{l_2}(x) T_{l_3}(x)}{\sqrt{Var(T_{l_1}(x)) Var(T_{l_2}(x)) Var(T_{l_3}(x))}} dx \qquad (12)$$

where

$$T_l(x) := L_l * T = \int_{S^2} T(y) \sum_m Y_{lm}(x) Y_{lm}^*(y) dy = \sum_m a_{lm} Y_{lm}(x),$$

i.e. we focus on the Fourier projections of the random fields at multipoles $(l_1, l_2, l_3)$. Both versions of the bispectrum have been shown to be extremely powerful against non-Gaussian alternatives, indeed there is now a widespread consensus that they make up the most efficient statistical procedures to search for non-Gaussianity, at least in the presence of fully observed spherical maps (see for instance [11, 21, 28]).

On the other hand, it is also well-known that the performance of Fourier methods in general, and the bispectrum in particular, deteriorates quite clearly in the presence of missing observations/partially observed maps ([20, 22, 11]). A natural idea is then to explore the localization properties of the needlets in harmonic domain, together with their robustness in the presence of gaps, in order to devise a statistic which might mimick the positive features of the bispectrum, at the same time coping with the difficulties brought in by missing observations.

To this aim, we shall consider the (normalized) *needlets bispectrum*

$$I_{j_1 j_2 j_3} := \sum_{k_3} \widehat{\beta}_{j_1 k_1} \widehat{\beta}_{j_2 k_2} \widehat{\beta}_{j_3 k_3} \delta_{j_1 j_2 j_3}(k_1, k_2, k_3) h_{j_1 j_2 j_3}(k_1, k_2, k_3),\ j_1 \le j_2 \le j_3,$$
$$(13)$$

where

$$\begin{aligned}
\delta_{j_1 j_2 j_3}(k_1, k_2, k_3) &= \mathbb{I}(\xi_{j_3 k_3} \in \mathcal{V}_{j_2 k_2}) \mathbb{I}(\xi_{j_2 k_2} \in \mathcal{V}_{j_1 k_1}), \\
h_{j_1 j_2 j_3}(k_1, k_2, k_3) &= \begin{cases} B^{j_2 - j_1} \frac{\sqrt{\lambda_{j_1 k_1}}}{\#\{k_2 : k_2 \in \mathcal{V}_{j_1 k_1} \cap \mathcal{X}_{j_2}\}}, & j_1 < j_2 = j_3, \\ \sqrt{\lambda_{j_3 k_3}}, & otherwise \end{cases},
\end{aligned}$$

where $\mathbb{I}(.)$ denotes the indicator function and $\mathcal{V}_{jk}$ is the Voronoi Cell that corresponds to $\xi_{jk}$; note that for $j_2 = j_3$ we have $h_{j_1 j_2 k_2} \approx \sqrt{\lambda_{j_2 k_2}}$. It is immediate to see that (13) can be seen as a natural development of (12), where the convolution with the orthonormal projector operator $L_l$ is replaced by the (discretized) convolution with the frame operator projection $\sqrt{\lambda_{jk}} \sum_l b(l/B^j) L_l$. Of course, in practice (12) is unfeasible and requires discretization to be implemented. In words, we are considering a version of the bispectrum where the exact identification of the multipoles is blurred by a form of suitable smoothing, with the purpose of a better robustness against missing observations.

The summation convention in (13) needs some further discussion. In practice, for applications the Voronoi tessellation is chosen in such a way to be nested across different scales (this is the case, for instance, for HEALPix [18], which is the standard package for CMB applications). Under such circumstances, our procedure can be described more explicitly as follows: we take

$$I_{j_1 j_2 j_3} := \sum_{k_3} \widehat{\beta}_{j_1 k_1} \widehat{\beta}_{j_2 k_2} \widehat{\beta}_{j_3 k_3} h_{j_1 j_2 j_3}(k_1, k_2, k_3), \tag{14}$$

where $k_2 = k_2(k_3)$ is the (unique) value of $k_2$ such that $\xi_{j_3 k_3} \in \mathcal{V}_{j_2 k_2}$, and $k_1 = k_1(k_3)$ is the (unique) value of $k_1$ such that $\xi_{j_2 k_2}, \xi_{j_3 k_3} \in \mathcal{V}_{j_1 k_1}$. In other words, the "finest grid" $\mathcal{X}_{j_3}$ is the one which leads the summation, whereas smaller frequency terms are identified with those *centres* whose corresponding Voronoi cells include the points being summed. Note, however, that in the general non-nested case the centre of $\mathcal{V}_{j_1 k_1}$ needs not belong to $\mathcal{X}_{j_3}$. In the sequel, for notational simplicity we write $k_1, k_2$ rather than $k_1(k_3), k_2(k_3)$, when no ambiguity is possible.

To investigate the asymptotic behaviour of the needlets bispectrum, we shall make an extensive use of the Diagram Formula which was introduced in the previous section. A crucial point, of course is the determination of the frequencies where the needlets bispectrum is evaluated. We distinguish three cases, namely:

- **Case 1)** $j_1 + 1 < j_2 < j_3 - 1$;
- **Case 2)** $j_1 + 1 < j_2 = j_3$, or $j_1 = j_2 < j_3 - 1$;
- **Case 3)** $j_1 = j_2 = j_3$.

Case 1 corresponds to the situation where all three frequencies are different. Case 2) is basically the "squeezed" or collapsed configuration which is considered by [3, 23], and many other cosmological references; in words, one frequency is (much) smaller than the other two. It has been widely argued in the physical literature that this configuration corresponds to the highest power region for so-called local models of non-Gaussianity. Case 3 corresponds to so-called equilateral configurations; this case, however, can be largely investigated by means of results in [6] and we report it only for completeness, omitting many details in the proof. It should be noted that for case 1) and 2) we focus on frequencies that are at least two steps apart, in order to exploit the semiorthogonality properties of the needlets systems. Relaxing this assumption implies no new ideas and would only make the paper notationally more complicated.

In each of the three cases we have trivially $\mathbb{E}I_{j_1 j_2 j_3} = 0$. We now focus on the asymptotic behaviour; here, asymptotics must be understood in the high resolution sense, i.e. we focus on a single realization of an isotropic random field, and we investigate the behaviour of our statistics at higher and higher bands. The first task is to ensure the statistics are non-degenerate and do not exhibit an explosive behaviour; this is the aim of the next Lemma. While the bound from above is quite straightforward, the lower bound is much more complicated and settles a question which was left open in [5], where the lower limit was simply assumed to be strictly larger than zero even for the simple case where $j_1 = j_2 = j_3$. As it is evident from the proof, the integer $K$ depends on the choice of cubature points and of kernel function $b(.)$; more explicit expressions can be found below. Note also that unless the three bands are equal, condition b) cannot be satisfied for $B = 2$; indeed in CMB applications values of order $B \simeq 1.2, 1.3$ are likely to be favoured.

**Lemma 1.** *Under Condition A, as $j_3 \to \infty$,*
*a) For all $j_1 \leq j_2 \leq j_3$, $\mathbb{E}I_{j_1 j_2 j_3}^2 = O(1)$.*
*Also, there exist a positive integer $K$ such that*
*b) For $\{j_1 = j_2 = j_3\}$, or $\{j_1 + K < j_2,\ B^{j_1} + B^{j_2} \geq B^{j_3}\}$, $\{\mathbb{E}I_{j_1 j_2 j_3}^2\}^{-1} = O(1)$.*
*The bounds are uniform with respect to $j_1, j_2$.*

*Proof.* **(a)** In the sequel, we shall use the fact that the set of the cubature points of polynomial spaces with degree less than $B^j$ are a $\kappa B^{-j}$−net; we also define $\rho := \max_{j,k} \{B^{-2j}\lambda_{j,k}\}$. The proof of all three cases are similar; we shall focus on case 2) $j_1 + 1 < j_2 = j_3$.

$$\mathbb{E}I_{j_1 j_2 j_2}^2$$
$$= \sum_{k_2, k_2' \in \mathcal{X}_{j_2}} \frac{\mathbb{E}\beta_{j_1 k_1}\beta_{j_1 k_1'} \mathbb{E}(\beta_{j_2 k_2}^2 \beta_{j_2 k_2'}^2)}{(\sigma_{j_1}\sigma_{j_2}\sigma_{j_2})^2} h_{j_1 j_2 j_3}(k_1, k_2, k_3) h_{j_1 j_2 j_3}(k_1', k_2', k_3')$$
$$= B^{2j_2 - 2j_1} \sum_{k_1, k_1' \in \mathcal{X}_{j_1}} \sum_{\substack{k_2 \in \mathcal{V}_{j_1 k_1} \cap \mathcal{X}_{j_2}, \\ k_2' \in \mathcal{V}_{j_1 k_1'} \cap \mathcal{X}_{j_2}}} \left( \mathbb{E}\widehat{\beta}_{j_1 k_1}\widehat{\beta}_{j_1 k_1'} + 2\mathbb{E}\widehat{\beta}_{j_1 k_1}\widehat{\beta}_{j_1 k_1'} (\mathbb{E}\widehat{\beta}_{j_2 k_2}\widehat{\beta}_{j_2 k_2'})^2 \right)$$
$$\times \frac{\sqrt{\lambda_{j_1 k_1}}}{\#\{k_2, k_2 \in \mathcal{V}_{j_1 k_1} \cap \mathcal{X}_{j_2}\}} \frac{\sqrt{\lambda_{j_1 k_1'}}}{\#\{k_2', k_2' \in \mathcal{V}_{j_1 k_1'} \cap \mathcal{X}_{j_2}\}}$$
$$= B^{2j_2 - 2j_1} \sum_{k_2 k_2'} \frac{\mathbb{E}\beta_{j_1 k_1}\beta_{j_1 k_1'} (\mathbb{E}\beta_{j_2 k_2}\beta_{j_2 k_2'})^2}{(\sigma_{j_1}\sigma_{j_2}\sigma_{j_2})^2} \frac{\sqrt{\lambda_{j_1 k_1}}}{\#\{k_2, k_2 \in \mathcal{V}_{j_1 k_1} \cap \mathcal{X}_{j_2}\}}$$
$$\times \frac{\sqrt{\lambda_{j_1 k_1'}}}{\#\{k_2', k_2' \in \mathcal{V}_{j_1 k_1'} \cap \mathcal{X}_{j_2}\}}.$$

Since $\lambda_{jk} \simeq B^{-2j}$ for every $\xi_{jk}$, and $\#\{k_2, k_2 \in \mathcal{V}_{j_1 k_1} \cap \mathcal{X}_{j_2}\} \simeq B^{2j_2 - 2j_1}$, the sum can be readily bounded by

$$\frac{C}{B^{2j_2}} \sum_{k_2 k_2'} \mathbb{E}\widehat{\beta}_{j_1 k_1}\widehat{\beta}_{j_1 k_1'} \left(\mathbb{E}\widehat{\beta}_{j_2 k_2}\widehat{\beta}_{j_2 k_2'}\right)^2 \quad \leq \quad \frac{C'}{B^{2j_2}} \sum_{k_2 k_2'} \left(\mathbb{E}\widehat{\beta}_{j_2 k_2}\widehat{\beta}_{j_2 k_2'}\right)^2$$
$$= \quad O(1) \text{ as } j_2 \to \infty,$$

in view of (10). This completes the proof of part a).

**(b)** The proof of the lower bound on the variance is considerably more delicate. We recall the correlation of the needlets coefficient is provided by

$$\mathbb{E}\widehat{\beta}_{jk}\widehat{\beta}_{jk'} \approx \frac{\sum_{l = B^{j-1}}^{B^{j+1}} b^2\left(\frac{l}{B^j}\right) C_l P_l(\cos\theta)}{\sum_{l = B^{j-1}}^{B^{j+1}} b^2\left(\frac{l}{B^j}\right) C_l \frac{2l+1}{4\pi}}$$

where $\theta = \arccos \langle \xi_{jk}, \xi_{jk'} \rangle$. The idea of our argument is to replace the needlets coefficients in the coarsest grid $\mathcal{X}_{j_1}$ by coefficients with the same resolution but evaluated over the pixels $\mathcal{X}_{j_2} \cap V_{j_1}(\xi_{j_1k_1})$. This will allow us to circumvent the fact that the cubature weights at the smaller frequencies stay constant over a Voronoi cell, whereas those corresponding to the highest $j$'s may vary. We shall hence consider

$$\beta_{jk}^*(\xi_{j'k'}) = \sqrt{\lambda_{jk}} \sum_{l=B^{j-1}}^{B^{j+1}} \sum_{m=-l}^{l} b\left(\frac{l}{B^j}\right) a_{lm} Y_{lm}(\xi_{j'k'}).$$

For fixed $(j', k')$, $j' < j$, $\beta_{jk}^*(\xi_{j'k'})$ varies over the pixels in $\mathcal{X}_j \cap V_{j'}(\xi_{j'k'})$. Let us now recall a few properties of Legendre polynomials, that we shall use extensively soon (see [36] for details); we have

$$\sup_{\theta \in [0,\pi]} P_l(\cos\theta) = P_l(\cos 0) = 1, \text{ and } \sup_{\theta \in [0,\pi]} \left| \frac{d}{d\theta} P_l(\cos\theta) \right| \leq 3l.$$

As a consequence, for any positive $\epsilon < 1$, there exists a $\delta > 0$, s.t. if $0 \leq \theta < \delta \leq \epsilon/(3l)$, then,

$$|P_l(\cos(\theta_0 + \theta)) - P_l(\cos\theta_0)| \leq 3l\theta \leq \epsilon,$$

because

$$0 \leq \cos(\theta_0 + \theta) - \cos\theta_0 = 2\sin\frac{2\theta_0 + \theta}{2}\sin\frac{\theta}{2} \leq \theta.$$

Now fix an integer $K$ such that $K \geq \log_B 6\kappa/\epsilon$; if we let $\{\xi_{j,k}\}_{j,k}$ be the cubature points of polynomial space with degree less than $B^{j+K+1}$, (note that we can assume all of these sets are $\kappa B^{-j}$-nets), then for any $j_2 > j_1, \xi_{j_2k_2} \in \mathcal{V}_{j_1k_1}, \xi_{j_2k_2'} \in \mathcal{V}_{j_1k_1'}$,

$$\begin{aligned}
\left| \langle \xi_{j_1k_1}, \xi_{j_1k_1'} \rangle - \langle \xi_{j_2k_2}, \xi_{j_2k_2'} \rangle \right| &= \left| d(\xi_{j_1k_1}, \xi_{j_1k_1'}) - d\left(\xi_{j_2k_2}, \xi_{j_2k_2'}\right) \right| \\
&\leq 2\kappa B^{-(j_1+K+1)} \leq B^{-(j_1+1)}\epsilon/3.
\end{aligned}$$

It follows that

$$\begin{aligned}
&\left| \mathbb{E}\beta_{j_1k_1}\beta_{j_1k_1'} - \mathbb{E}[\beta_{j_1k_1}^*(\xi_{j_2k_2})\beta_{j_1k_1'}^*(\xi_{j_2k_2'})] \right| \\
&\leq \sqrt{\lambda_{jk}\lambda_{jk'}} \sum_{l=B^{j_1-1}}^{B^{j_1+1}} b^2\left(\frac{l}{B^j}\right) \\
&\quad \times C_l \left| \sum_{m=-l}^{l} Y_{lm}(\xi_{j_1k_1})\overline{Y_{lm}(\xi_{j_1k_1'})} - \sum_{m=-l}^{l} Y_{lm}(\xi_{j_2k_2})\overline{Y_{lm}(\xi_{j_2k_2'})} \right| \\
&= \sqrt{\lambda_{jk}\lambda_{jk'}} \sum_{l=B^{j_1-1}}^{B^{j_1+1}} b^2\left(\frac{l}{B^j}\right) C_l \frac{2l+1}{4\pi} \left| P_l\left(\left\langle \xi_{j_1k_1}, \xi_{j_1k_1'} \right\rangle\right) - P_l\left(\left\langle \xi_{j_2k_2}, \xi_{j_2k_2'} \right\rangle\right) \right| \\
&\leq \epsilon\sqrt{\lambda_{jk}\lambda_{jk'}} \sum_{l=B^{j_1-1}}^{B^{j_1+1}} b^2\left(\frac{l}{B^j}\right) C_l \frac{2l+1}{4\pi} \leq \epsilon C \sigma_{j_1}^2.
\end{aligned} \tag{15}$$

where $C$ is a constant, $C = C(\kappa)$ (see [29]). We are now in the position to establish our lower bound. By simple algebraic manipulations, we have

$$
\mathbb{E} I_{j_1 j_2 j_3}^2
$$
$$
\approx \sum_{k_3, k_3' \in \mathcal{X}_{j_3}} \frac{\mathbb{E}\beta_{j_1 k_1}^* (\xi_{j_3 k_3}) \beta_{j_1 k_1'}^* (\xi_{j_3 k_3'}) \mathbb{E}\beta_{j_2 k_2}^* (\xi_{j_3 k_3}) \beta_{j_2 k_2'}^* (\xi_{j_3 k_3'}) \mathbb{E}\beta_{j_3 k_3} \beta_{j_3 k_3'}}{(\sigma_{j_1} \sigma_{j_2} \sigma_{j_3})^2}
$$
$$
\times \sqrt{\lambda_{j_3 k_3}} \sqrt{\lambda_{j_3 k_3'}}
$$
$$
+ \sum_{k_3, k_3' \in \mathcal{X}_{j_3}} \frac{\left( \mathbb{E}\beta_{j_1 k_1} \beta_{j_1 k_1'} - \mathbb{E}\beta_{j_1 k_1}^* (\xi_{j_3 k_3}) \beta_{j_1 k_1'}^* (\xi_{j_3 k_3'}) \right)}{\sigma_{j_1}^2}
$$
$$
\times \mathbb{E}\widehat{\beta}_{j_2 k_2} \widehat{\beta}_{j_2 k_2'} \mathbb{E}\widehat{\beta}_{j_3 k_3} \widehat{\beta}_{j_3 k_3'} \sqrt{\lambda_{j_3 k_3}} \sqrt{\lambda_{j_3 k_3'}}
$$
$$
+ \sum_{k_3, k_3' \in \mathcal{X}_{j_3}} \frac{\mathbb{E}\beta_{j_1 k_1}^* (\xi_{j_3 k_3}) \beta_{j_1 k_1'}^* (\xi_{j_3 k_3'}) \left( \mathbb{E}\beta_{j_2 k_2} \beta_{j_2 k_2'} - \mathbb{E}\beta_{j_2 k_2}^* (\xi_{j_3 k_3}) \beta_{j_2 k_2'}^* (\xi_{j_3 k_3'}) \right)}{(\sigma_{j_1} \sigma_{j_2})^2}
$$
$$
\times \mathbb{E}\widehat{\beta}_{j_3 k_3} \widehat{\beta}_{j_3 k_3'} \sqrt{\lambda_{j_3 k_3}} \sqrt{\lambda_{j_3 k_3'}}.
$$

In view of (15), (10) and standard manipulations, the second and third components can be made arbitrarily small. To bound the first component, we recall first two facts involving spherical harmonics (see again ([36, 23] for details), namely the so-called Gaunt integral

$$
\int_{S^2} Y_{l_1 m_1}(x) Y_{l_2 m_2}(x) Y_{l_3 m_3}(x) dx
$$
$$
= \sqrt{\frac{(2l_1+1)(2l_2+1)(2l_3+1)}{4\pi}} \begin{pmatrix} l_1 & l_2 & l_3 \\ m_1 & m_2 & m_3 \end{pmatrix} \begin{pmatrix} l_1 & l_2 & l_3 \\ 0 & 0 & 0 \end{pmatrix},
$$

where for the Wigner's 3j coefficients introduced on the right-hand side we recall the properties

$$
\sum_{m_1 m_2 m_3} \begin{pmatrix} l_1 & l_2 & l_3 \\ m_1 & m_2 & m_3 \end{pmatrix}^2 \equiv 1, \tag{16}
$$

$$
\begin{pmatrix} l_1 & l_2 & l_3 \\ 0 & 0 & 0 \end{pmatrix}^2 \geq \frac{C}{l_3^2} \mathbb{I}\{l_1 + l_2 \geq l_3, l_1 + l_2 + l_3 = even\}. \tag{17}
$$

Using the fact that $\{\lambda_{jk}\}$ are the cubature weights corresponding to the space $\mathcal{K}_{3B^{j+1}}$, easy manipulations yield

$$
\frac{1}{\sigma_{j_1}^2 \sigma_{j_2}^2 \sigma_{j_3}^2} B^{2j_3} \sum_{l_1 l_2 l_3} \sum_{m_1 m_2 m_3} b^2 \left( \frac{l_1}{B^{j_1}} \right) b^2 \left( \frac{l_2}{B^{j_2}} \right) b^2 \left( \frac{l_3}{B^{j_3}} \right) C_{l_1} C_{l_2} C_{l_3}
$$
$$
\times \sum_{k_3} Y_{l_1 m_1}(\xi_{j_3 k_3}) Y_{l_2 m_2}(\xi_{j_3 k_3}) Y_{l_3 m_3}(\xi_{j_3 k_3}) \lambda_{j_3 k_3}
$$
$$
\times \sum_{k_3'} Y_{l_1 m_1}(\xi_{j_3 k_3'}) Y_{l_2 m_2}(\xi_{j_3 k_3'}) Y_{l_3 m_3}(\xi_{j_3 k_3'}) \lambda_{j_3 k_3'}
$$

$$= \frac{1}{\sigma_{j_1}^2 \sigma_{j_2}^2 \sigma_{j_3}^2} B^{2j_3} \sum_{l_1 l_2 l_3} \sum_{m_1 m_2 m_3} b^2\left(\frac{l_1}{B^{j_1}}\right) b^2\left(\frac{l_2}{B^{j_2}}\right) b^2\left(\frac{l_3}{B^{j_3}}\right) C_{l_1} C_{l_2} C_{l_3}$$

$$\times \left(\int_{S^2} Y_{l_1 m_1}(x) Y_{l_2 m_2}(x) Y_{l_3 m_3}(x) dx\right)^2$$

$$= \frac{1}{\sigma_{j_1}^2 \sigma_{j_2}^2 \sigma_{j_3}^2} B^{2j_3} \sum_{l_1 l_2 l_3} b^2\left(\frac{l_1}{B^{j_1}}\right) b^2\left(\frac{l_2}{B^{j_2}}\right) b^2\left(\frac{l_3}{B^{j_3}}\right)$$

$$\times C_{l_1} C_{l_2} C_{l_3} \frac{(2l_1+1)(2l_2+1)(2l_3+1)}{4\pi}$$

$$\times \sum_{m_1 m_2 m_3} \left(\begin{array}{ccc} l_1 & l_2 & l_3 \\ m_1 & m_2 & m_3 \end{array}\right)^2 \left(\begin{array}{ccc} l_1 & l_2 & l_3 \\ 0 & 0 & 0 \end{array}\right)^2$$

and using $\sigma_j^2 \approx B^{2j} C_{B^j}$, $B^{j_i} \approx l_i, i = 1, 2, 3$ we get

$$\frac{1}{\sigma_{j_1}^2 \sigma_{j_2}^2 \sigma_{j_3}^2} B^{2j_3} \sum_{l_1 l_2 l_3} b^2\left(\frac{l_1}{B^{j_1}}\right) b^2\left(\frac{l_2}{B^{j_2}}\right) b^2\left(\frac{l_3}{B^{j_3}}\right)$$

$$\times C_{l_1} C_{l_2} C_{l_3} \frac{(2l_1+1)(2l_2+1)(2l_3+1)}{4\pi} \left(\begin{array}{ccc} l_1 & l_2 & l_3 \\ 0 & 0 & 0 \end{array}\right)^2$$

$$\approx B^{-j_1-j_2-j_3} \sum_{\substack{l_1+l_2 \geq l_3 \\ l_1+l_2+l_3=even}} b^2\left(\frac{l_1}{B^{j_1}}\right) b^2\left(\frac{l_2}{B^{j_2}}\right) b^2\left(\frac{l_3}{B^{j_3}}\right) > c > 0.$$

In the previous argument, we have taken $c := \inf_{j_1 j_2 j_3} c_{j_1 j_2 j_3}$, where

$$c_{j_1 j_2 j_3} := B^{-j_1-j_2-j_3} \sum_{\substack{l_1+l_2 \geq l_3 \\ l_1+l_2+l_3=even}} b^2\left(\frac{l_1}{B^{j_1}}\right) b^2\left(\frac{l_2}{B^{j_2}}\right) b^2\left(\frac{l_3}{B^{j_3}}\right).$$

It is simple to show that $c > 0$. Indeed, under (4), we have

$$\begin{aligned} c_{j_1 j_2 j_3} &\geq \frac{1}{2} \int_{\substack{\{u_1+u_2 \geq u_3\}, \\ [0,2]^3 \cap \{u_1 \leq u_2 \leq u_3\}}} b^2(u_1) b^2(u_2) b^2(u_3) du_1 du_2 du_3 \\ &\geq \frac{1}{2} \int_{1-\delta_0}^{1+\delta_0} du_1 \int_{u_1}^{1+\delta_0} du_2 \int_{1+\delta_0}^{2-2\delta_0} du_3 \{b^2(u_1) b^2(u_2) b^2(u_3)\} \\ &\geq \frac{c_0}{8} \delta_0^3, \end{aligned}$$

where $\delta_0 \leq 1/4$, and $b(x) \geq 1/2$ for any $x \in [1-\delta_0, 1+\delta_0]$, $c_0 = \inf_{x \in [1+\delta_0, 2-2\delta_0]} b^2(x)$.

The same argument as before could be used to establish lower bounds when $j_1 = j_2 < j_3$ or $j_1 = j_2 = j_3$. To conclude our proof, we consider the case when $j_1 < j_2 = j_3$. We obtain

$$B^{2j_2-2j_1} \sum_{k_2 k_2'} \mathbb{E}\widehat{\beta}_{j_1 k_1} \widehat{\beta}_{j_1 k_1'} \left( \mathbb{E}\widehat{\beta}_{j_2 k_2} \widehat{\beta}_{j_2 k_2'} \right)^2$$

$$\times \frac{\sqrt{\lambda_{j_1 k_1}}}{\#\{k_2, k_2 \in \mathcal{V}_{j_1 k_1} \cap \mathcal{X}_{j_2}\}} \frac{\sqrt{\lambda_{j_1 k_1'}}}{\#\{k_2', k_2' \in \mathcal{V}_{j_1 k_1'} \cap \mathcal{X}_{j_2}\}}$$

$$\geq C\left(\kappa\right) \frac{B^{2j_2-2j_1}}{B^{4j_2-4j_1}} \sum_{k_1} \sum_{k_2 \in \mathcal{V}_{j_1 k_1} \cap \mathcal{X}_{j_2}} \left( \mathbb{E}\widehat{\beta}_{j_2 k_2}^2 \right)^2 \lambda_{j_1 k_1}$$

$$+ C'\left(\kappa, \rho\right) \frac{B^{2j_2-2j_1}}{B^{4j_2-4j_1}} \sum_{k_1 \neq k_1'} \sum_{\substack{k_2 \in \mathcal{V}_{j_1 k_1} \cap \mathcal{X}_{j_2}, \\ k_2' \in \mathcal{V}_{j_1 k_1'} \cap \mathcal{X}_{j_2}}} \left( \mathbb{E}\widehat{\beta}_{j_2 k_2} \widehat{\beta}_{j_2 k_2'} \right)^2,$$

the two summands corresponding to the cases where $k_2, k_2'$ belong to the same or to different Voronoi cells, respectively. The first part is equal to

$$C\left(\kappa\right) B^{2j_1-2j_2} \sum_{k_1} \sum_{k_2 \in \mathcal{V}_{j_1 k_1} \cap \mathcal{X}_{j_2}} \lambda_{j_1 k_1} \geq C\left(\kappa\right) \sum_{k_1} \lambda_{j_1 k_1} = 4\pi C\left(\kappa\right),$$

while the second part is equal to

$$C'\left(\kappa, \rho\right) B^{-2j_2} \left\{ \sum_{\substack{k_1 \neq k_1'; \\ d(\mathcal{V}_{j_1 k_1}, \mathcal{V}_{j_1 k_1'})=0}} \sum_{\substack{k_2 \in \mathcal{V}_{j_1 k_1} \cap \mathcal{X}_{j_2}, \\ k_2' \in \mathcal{V}_{j_1 k_1'} \cap \mathcal{X}_{j_2}}} + \sum_{\substack{k_1 \neq k_1'; \\ d(\mathcal{V}_{j_1 k_1}, \mathcal{V}_{j_1 k_1'})>0}} \sum_{\substack{k_2 \in \mathcal{V}_{j_1 k_1} \cap \mathcal{X}_{j_2}, \\ k_2' \in \mathcal{V}_{j_1 k_1'} \cap \mathcal{X}_{j_2}}} \right\}$$

$$\times \left( \mathbb{E}\widehat{\beta}_{j_2 k_2} \widehat{\beta}_{j_2 k_2'} \right)^2$$

$$\leq C'\left(\kappa, \rho\right) \frac{6\pi}{B^{2j_2}} \sum_{k_1} \sum_{\substack{k_2 \in \mathcal{V}_{j_1 k_1} \cap \mathcal{X}_{j_2}, \\ k_2' \in \mathcal{V}_{j_1 k_1'} \cap \mathcal{X}_{j_2}}} \left( \mathbb{E}\widehat{\beta}_{j_2 k_2} \widehat{\beta}_{j_2 k_2'} \right)^2$$

$$+ \frac{C'\left(\kappa, \rho\right)}{B^{2j_2}} \sum_{k_1} \sum_{k_2 \in \mathcal{V}_{j_1 k_1} \cap \mathcal{X}_{j_2}} \left\{ \sum_{k_1'} \sum_{\substack{k_2' \in \mathcal{V}_{j_1 k_1'} \cap \mathcal{X}_{j_2}; \\ d(\mathcal{V}_{j_1 k_1}, \mathcal{V}_{j_1 k_1'})>0}} \left( \mathbb{E}\widehat{\beta}_{j_2 k_2} \widehat{\beta}_{j_2 k_2'} \right)^2 \right\}$$

$$\leq C'\left(\kappa, \rho\right) \frac{6\pi}{B^{2j_2}} \sum_{k_1} B^{j_2-j_1}$$

$$+ \frac{C'\left(\kappa, \rho\right)}{B^{2j_2}} \sum_{k_2} \sum_{k_2'} \frac{C_M}{\left(1 + B^{j_2}(B^{-j_1} + d(k_2', \mathcal{V}_{j_1 k_1 (k_2)}))\right)^M}$$

$$\leq C'\left(\kappa, \rho\right) B^{-(j_2-j_1)} + C'\left(\kappa, \rho, C_M\right) \int_{B^{-j_1}}^{\infty} \frac{B^{2j_2} \sin\theta}{(1 + B^{j_2}\theta)^M} d\theta$$

$$= C\left(\kappa, \rho, C_M\right) B^{-(j_2-j_1)})$$

By taking $j_2 - j_1 \geq K \geq \max\left\{ [\log_B \{2\kappa/\epsilon\}], \frac{1}{2} c_1/c_2 \right\}$, where $c_1 = C(\kappa), c_2 = C\left(\kappa, \rho, C_M\right)$, and $\epsilon \leq \frac{1}{2} c/C_M$, the argument is completed. The proof for case (3) is similar (actually slightly easier). $\square$

The following weak convergence theorem is the main result of this Section. We stress that the statement could be easily extended to a multivariate Central Limit Theorem; however, because this extension would not entail any substantial novelty, at the same time making the notation much more cumbersome, we prefer to stick to the univariate case.

**Theorem 2.** *Let $T(x)$ be a zero-mean, mean square continuous and isotropic Gaussian random field, with angular power spectrum that satisfies Condition A. As $j_1 \to \infty$, we have*

$$\frac{I_{j_1 j_2 j_3}}{\sqrt{\mathbb{E} I_{j_1 j_2 j_3}^2}} \to_d N(0,1).$$

*Proof.* In view of the results in [32], to establish a Central Limit Theorem for a multilinear form in Gaussian random variables, it is enough to investigate the asymptotic behaviour of fourth order moments (or equivalently cumulants), see also [14]. Our aim is then to show that, as $j_1 \to \infty$,

$$\mathbb{E} I_{j_1 j_2 j_3}^4 = 3(\mathbb{E} I_{j_1 j_2 j_3}^2)^2 + O(B^{-j_1/2}).$$

For notational simplicity, we shall write

$$\rho_j(k', k) := \mathbb{E} \widehat{\beta}_{j,k'} \widehat{\beta}_{j,k}.$$

By the diagram formula we have, for all index sets $I$:

$$\mathbb{E} \left\{ \prod_{i \in I} \prod_{l=1}^3 \widehat{\beta}_{j_l k_l^i} \right\} = \sum_{\gamma \in \Gamma(I,3)} \prod_{\{(i,l)(i',l')\} \in \gamma} \delta_{j_l}^{j_{l'}} \rho_{j_l}(k_l^{(i)}, k_{l'}^{(i')}).$$

Similarly to [23], we define

$$\begin{aligned}
\rho(\gamma; j_1, j_2, j_3) &= \sum_{k_3^i \in \mathcal{X}_{j_3}} \prod_{\{(i,l)(i',l')\} \in \gamma} \delta_{j_l}^{j_{l'}} \rho_{j_l}(k_l^{(i)}, k_{l'}^{(i')}) \\
&\quad \times \prod_{i \in I} \delta_{j_1 j_2 j_3}(k_1^{(i)}, k_2^{(i)}, k_3^{(i)}) h_{j_1 j_2 j_3}(k_1^{(i)}, k_2^{(i)}, k_3^{(i)})
\end{aligned}$$

so that

$$\begin{aligned}
\mathbb{E} I_{j_1 j_2 j_3}^4 &= \mathbb{E} \left\{ \sum_{k_3} \widehat{\beta}_{j_1 k_1} \widehat{\beta}_{j_2 k_2} \widehat{\beta}_{j_3 k_3} \delta_{j_1 j_2 j_3}(k_1, k_2, k_3) h_{j_1 j_2 j_3}(k_1, k_2, k_3) \right\}^4 \\
&= \sum_{k_3^{(1)}} \cdots \sum_{k_3^{(4)}} \mathbb{E} \left\{ \prod_{i=1}^4 \prod_{l=1}^3 \widehat{\beta}_{j_l k_l^i} \delta_{j_1 j_2 j_3}(k_1^{(i)}, k_2^{(i)}, k_3^{(i)}) h_{j_1 j_2 j_3}(k_1^{(i)}, k_2^{(i)}, k_3^{(i)}) \right\} \\
&= \left\{ \sum_{\gamma \in \Gamma_C(4,3)} + \sum_{\gamma \in \Gamma_{\overline{C}}(4,3)} \right\} \rho(\gamma; j_1, j_2, j_3)
\end{aligned}$$

where $\Gamma_C$ is the set of all connected diagrams. To conclude our argument, we only need to show that

$$\sum_{\gamma \in \Gamma_C(4,3)} \rho(\gamma; j_1, j_2, j_3) = O(B^{-j_1/2}), \text{ as } j_1 \to \infty, \tag{18}$$

and

$$\sum_{\gamma \in \Gamma_{\overline{CF}}(4,3)} \rho(\gamma; j_1, j_2, j_3) = \sum_{\gamma \in \Gamma_P(4,3)} \rho(\gamma; j_1, j_2, j_3) = 3(\mathbb{E}I_{j_1 j_2 j_3}^2)^2. \qquad (19)$$

(19) is an immediate consequence of the definition of $\mathbb{E}I_{j_1 j_2 j_3}^2$ and trivial combinatorial manipulations(see [34]). The result in (18) is proved by splitting connected diagrams into those with or without flat edges. Diagrams with flat edges are dealt with in Lemma 3, while those without are dealt with in Lemma 4. We stress that, on the contrary of what is often the case when the diagram formula is applied, terms with flat edges do not vanish, due to correlation among different locations in the spherical needlets coefficients. This completes the proof of the main result. □

**Lemma 3.** *For a connected diagram with flat edges, $\gamma \in \Gamma_{CF}(4,3)$, we have*

$$\rho(\gamma; j_1, j_2, j_3) = O(B^{-2j_1}), \text{ as } j_1 \to \infty.$$

*Proof.* We write $\{(k_b^{(a)}, j_b)\}_{a=1,\dots,4,b=1,2,3}$ for the elements in our diagram, $a$ and $b$ denoting the row and column indexes, respectively. We recall also that $k_2 = k_2(k_3)$, $k_1 = k_1(k_3)$, as explained earlier. For Case 1), i.e. $j_1 < j_2 - 1 < j_3 - 2$, since $E\{\beta_{j_3 k_3} \beta_{j_2 k_2}\} = 0$ for every $k_2 \in \mathcal{X}_{j_2}, k_3 \in \mathcal{X}_{j_3}$,it is easy to see that $\rho(\gamma; j_1, j_2, j_3) \equiv 0$. For Case 2), i.e. $j_1 + 1 < j_2 = j_3$, we assume (without loss of generality) that a flat edge is present in the first row of the diagram, i.e. we let $\{(k_2^{(1)}, j_2)(k_3^{(1)}, j_3)\} \in \gamma$. By the same argument as in (7) we obtain immediately

$$B^{4j_2 - 4j_1} \prod_{i \in \{2,\dots,4\}} \sum_{k_2^{(i)}, i \neq 1} \prod_{\substack{\{(i,l)(i',l')\} \in \gamma, \\ i,i' \neq 1}} \delta_{ji}^{ji'} \rho_{j_l}(k_l^{(i)}, k_{l'}^{(i')})$$

$$\times \sum_{k_2^{(1)}} \frac{\rho_{j_1}(k_1^{(1)}, k_1^{(2)}) \sqrt{\lambda_{j_1 k_1^{(1)}}}}{\#\{k_2^{(1)}, k_2^{(1)} \in \mathcal{V}_{j_1 k_1^{(1)}} \cap \mathcal{X}_{j_2}\}} = 0.$$

On the other hand, under $j_1 = j_2 < j_3 - 1$, again we assume a flat edge $\{(k_1^{(1)}, j_1)(k_1^{(1)}, j_2)\} \in \gamma$; by necessity, there should exist another flat edge in the graph, and w.l.o.g. we take it be in the fourth row, i.e. $\{(k_1^{(4)}, j_1)(k_2^{(4)}, j_2)\} \in \gamma$. Then we have $\xi_{j_1 k_1^{(i)}} = \xi_{j_2 k_2^{(i)}}$, and the resulting term can be bounded by

$$\frac{1}{B^{4j_3}} \prod_{i \in \{1,\dots,4\}} \sum_{k_3^{(i)}} \left| \rho_{j_3}(k_3^{(1)}, k_3^{(2)}) \rho_{j_3}(k_3^{(3)}, k_3^{(4)}) \right| \rho_{j_1}^2(k_1^{(2)}, k_1^{(3)})$$

$$= \frac{1}{B^{4j_3}} \sum_{k_3^{(2)}, k_3^{(3)}} \rho_{j_1}^2(k_1^{(2)}, k_1^{(3)}) \left\{ \sum_{k_3^{(1)}} \left| \rho_{j_3}(k_3^{(1)}, k_3^{(2)}) \right| \sum_{k_3^{(4)}} \left| \rho_{j_3}(k_3^{(3)}, k_3^{(4)}) \right| \right\}$$

$$\leq \frac{1}{B^{4j_3}} \sum_{k_3^{(2)}, k_3^{(3)}} \rho_{j_1}^2(k_1^{(2)}, k_1^{(3)}) C_M \times C_M$$

$$\leq \frac{C}{B^{4j_3}} \sum_{k_1^{(2)}, k_1^{(3)}} \rho_{j_2}^2(k_1^{(2)}, k_1^{(3)}) \left[ \max_{k_1^{(2)}} \# \left\{ \xi_{j_3 k_3} \in \mathcal{V}_{j_1}(k_1^{(2)}) \right\} \right]$$

$$\times \left[ \max_{k_1^{(3)}} \# \left\{ \xi_{j_3 k_3} \in \mathcal{V}_{j_1}(k_1^{(3)}) \right\} \right]$$

$$= \frac{CB^{2j_1}}{B^{4j_3}} B^{4j_3 - 4j_1} = O(B^{-2j_1}).$$

Finally for Case 3), the argument is analogous; more precisely, components with diagrams in $\Gamma_{CF}(4,3)$ can be bounded by

$$|\rho(\gamma; j, j, j)| = \frac{1}{B^{4j}} \prod_{i \in I} \sum_{k^{(i)}} \left| \rho_j(k^{(1)}, k^{(2)}) \rho_j(k^{(3)}, k^{(4)}) \right| \rho_j^2(k^{(2)}, k^{(3)})$$

$$= \frac{1}{B^{4j}} \sum_{k^{(2)}, k^{(3)}} \rho_j^2(k^{(2)}, k^{(3)}) \left\{ \sum_{k^{(1)}} \left| \rho_j(k^{(1)}, k^{(2)}) \right| \sum_{k^{(4)}} \left| \rho_j(k^{(3)}, k^{(4)}) \right| \right\}$$

$$\leq \frac{1}{B^{4j}} C_M \sum_{k^{i_2}, k^{i_3}} \rho_j^2(k^{(2)}, k^{(3)}) = O(B^{-2j}),$$

where in the third equation we used (9). Thus the proof is completed. □

**Lemma 4.** *For a connected diagram without a flat edge, $\gamma \in \Gamma_{C\overline{F}}(4,3)$, we have*

$$\rho(\gamma; j_1, j_2, j_3) = O(B^{-j_1/2}), \text{ as } j_1 \to \infty.$$

*Proof.* Connected diagrams without flat edges and with four nodes can be partitioned into two classes, i.e. so-called *cliques,* where each vertex is connected to all three others, and terms with loops of order 2. We focus on the former class; without loss of generality, we can express $\rho(\gamma; j_1, j_2, j_3)$ by

$$\prod_{i=1}^{4} \sum_{k_3^i} \rho_{j_1}(k_1^{(1)}, k_1^{(2)}) \rho_{j_1}(k_1^{(3)}, k_1^{(4)}) \rho_{j_2}(k_2^{(1)}, k_2^{(3)}) \rho_{j_2}(k_2^{(2)}, k_2^{(4)})$$

$$\times \rho_{j_3}(k_3^{(1)}, k_3^{(4)}) \rho_{j_3}(k_3^{(2)}, k_3^{(3)})$$

$$\times \prod_{i=1}^{4} \delta_{j_1 j_2 j_3}(k_1^{(i)}, k_2^{(i)}, k_3^{(i)}) h_{j_1 j_2 j_3}(k_1^{(i)}, k_2^{(i)}, k_3^{(i)}).$$

By means of (9) we readily obtain the bound

$$\frac{C}{B^{4j_3}} \prod_{i=1}^{4} \sum_{k_3^i} \frac{1}{(1 + B^{j_1} d(k_1^{(1)}, k_1^{(2)}))^M} \frac{1}{(1 + B^{j_1} d(k_1^{(3)}, k_1^{(4)}))^M}$$

$$\times \frac{1}{(1 + B^{j_2} d(k_2^{(1)}, k_2^{(3)}))^M} \frac{1}{(1 + B^{j_2} d(k_2^{(2)}, k_2^{(4)}))^M}$$

$$\times \frac{1}{(1 + B^{j_3} d(k_3^{(1)}, k_3^{(4)}))^M} \frac{1}{(1 + B^{j_3} d(k_3^{(2)}, k_3^{(3)}))^M}. \tag{20}$$

Consider first Case 2), with $j_1 + 1 < j_2 = j_3$. Using inequality (11), (20) can be replaced by

$$B^{4j_2-4j_1} \sum_{k_2^{(1)},\ldots,k_2^{(4)} \in \mathcal{X}_{j_2}} \frac{1}{(1+B^{j_2}d(k_2^{(1)},k_2^{(3)}))^M} \frac{1}{(1+B^{j_2}d(k_2^{(2)},k_2^{(4)}))^M}$$

$$\times \frac{1}{(1+B^{j_2}d(k_2^{(1)},k_2^{(4)}))^M} \frac{1}{(1+B^{j_2}d(k_2^{(2)},k_2^{(3)}))^M}$$

$$\times \prod_{i=1}^{4} \frac{\sqrt{\lambda_{j_2 k_2^{(i)}}}}{\#\{k_2^{(i)}, k_2^{(i)} \in \mathcal{V}_{j_1 k_1^{(i)}} \cap \mathcal{X}_{j_2}\}}$$

$$\leq B^{4j_2} \frac{C}{B^{8j_2}} \sum_{k_2^{(3)},k_2^{(4)}} \left\{ \sum_{k_2^{(1)}} \frac{1}{(1+B^{j_2}d(k_2^{(1)},k_2^{(3)}))^M} \frac{1}{(1+B^{j_2}d(k_2^{(1)},k_2^{(4)}))^M} \right\}$$

$$\times \left\{ \sum_{k_2^{(2)}} \frac{1}{(1+B^{j_2}d(k_2^{(2)},k_2^{(4)}))^M} \frac{1}{(1+B^{j_2}d(k_2^{(2)},k_2^{(3)}))^M} \right\}$$

$$\leq \frac{C}{B^{4j_2}} \sum_{k_2^{(3)},k_2^{(4)}} \left\{ \frac{C_M}{(1+B^{j_2}d(k_2^{(3)},k_2^{(4)}))^{2M}} \right\}$$

$$\leq \frac{C}{B^{4j_2}} C_M C_{2M} B^{2j_2} = O(B^{-2j_2}) = o(B^{-2j_1}).$$

Likewise, for $j_1 = j_2 < j_3 - 1$, we have the bound

$$\frac{C}{B^{4j_3}} \prod_{i=1}^{4} \sum_{k_3^{(i)}} \frac{1}{(1+B^{j_1}d(k_1^{(1)},k_1^{(2)}))^M} \frac{1}{(1+B^{j_1}d(k_1^{(3)},k_1^{(4)}))^M}$$

$$\times \frac{1}{(1+B^{j_2}d(k_1^{(1)},k_1^{(3)}))^M} \frac{1}{(1+B^{j_2}d(k_1^{(2)},k_1^{(4)}))^M} \frac{1}{(1+B^{j_3}d(k_3^{(1)},k_3^{(4)}))^M}$$

$$\times \frac{1}{(1+B^{j_3}d(k_3^{(2)},k_3^{(3)}))^M}$$

$$\leq \frac{C}{B^{4j_3}} \prod_{i=1}^{4} \sum_{k_1^{(i)}} \frac{1}{(1+B^{j_1}d(k_1^{(1)},k_1^{(2)}))^M} \frac{1}{(1+B^{j_1}d(k_1^{(3)},k_1^{(4)}))^M}$$

$$\times \frac{1}{(1+B^{j_2}d(k_1^{(1)},k_1^{(3)}))^M} \frac{1}{(1+B^{j_2}d(k_1^{(2)},k_1^{(4)}))^M}$$

$$\times \sum_{k_3^1 \in \mathcal{V}_{j_1}(k_1^{(1)}) \cap \mathcal{X}_{j_3}} \frac{1}{(1+B^{j_3}d(k_3^{(1)},k_3^{(4)}))^M}$$

$$\times \sum_{k_3^2 \in \mathcal{V}_{j_1}(k_1^{(2)}) \cap \mathcal{X}_{j_3}} \frac{1}{(1+B^{j_3}d(k_3^{(2)},k_3^{(3)}))^M}$$

$$\leq \frac{CB^{4j_3-4j_1}}{B^{4j_3}} \sum_{k_1^{(2)},k_1^{(3)}} \left\{ \sum_{k_1^{(1)}} \frac{1}{(1+B^{j_1}d(k_1^{(1)},k_1^{(2)}))^M} \frac{1}{(1+B^{j_2}d(k_1^{(1)},k_1^{(3)}))^M} \right\}$$

$$\times \left\{ \sum_{k_1^{(4)}} \frac{1}{(1+B^{j_1}d(k_1^{(3)},k_1^{(4)}))^M} \frac{1}{(1+B^{j_2}d(k_1^{(2)},k_1^{(4)}))^M} \right\}$$

$$\leq \frac{CB^{4j_3-4j_1}}{B^{4j_3}} \sum_{k_1^{(2)},k_1^{(3)}} \frac{C_M^2}{(1+B^{j_2}d(k_1^{(2)},k_1^{(3)}))^{2M}} \leq CB^{2j_1-4j_1} = O(B^{-2j_1}).$$

This concludes the proof for Case 2); the proof for Case 3) could be implemented along identical lines.

The analysis for case 1) $j_1 + 1 < j_2 < j_3 - 1$, i.e. the case where all three frequencies differ, requires considerably more care. As before, let $\mathcal{V}_{j_u}(k_u^{(i)}, r)$ be the Voronoi cells associated to $k_u^{(i)}$, and we recall it satisfies $B(x_i, r/2) \subset \mathcal{V}_{j_u}(x_i) \subset B(x_i, r)$, $k_u^{(i)} \in \mathcal{X}_{j_u}$. Our idea is to partition (20) into four elements, as follows:

$$(20) \leq \frac{C}{B^{4j_3}} \prod_{i=1}^4 \sum_{k_3^i} \{ \mathbb{I}(d(k_1^{(1)}, k_1^{(4)} > 3r), d(k_1^{(2)}, k_1^{(3)} > 3r))$$

$$+ \mathbb{I}(d(k_1^{(1)}, k_1^{(4)} > 3r), d(k_1^{(2)}, k_1^{(3)} \leq 3r))$$

$$+ \mathbb{I}(d(k_1^{(1)}, k_1^{(4)} \leq 3r), d(k_1^{(2)}, k_1^{(3)} > 3r))$$

$$+ \mathbb{I}(d(k_1^{(1)}, k_1^{(4)} \leq 3r), d(k_1^{(2)}, k_1^{(3)} \leq 3r)) \}$$

$$\times \frac{1}{(1+B^{j_1}d(k_1^{(1)}, k_1^{(2)}))^M} \frac{1}{(1+B^{j_1}d(k_1^{(3)}, k_1^{(1)}))^M} \frac{1}{(1+B^{j_2}d(k_2^{(1)}, k_2^{(3)}))^M}$$

$$\times \frac{1}{(1+B^{j_2}d(k_2^{(2)}, k_2^{(4)}))^M} \frac{1}{(1+B^{j_3}d(k_3^{(1)}, k_3^{(4)}))^M} \frac{1}{(1+B^{j_3}d(k_3^{(2)}, k_3^{(3)}))^M}$$

The first three summands are easy to bound, indeed it is enough to notice that

$$\frac{C}{B^{4j_3}} \prod_{i=1}^4 \sum_{k_3^i} \left\{ \mathbb{I}\Big(d\left(k_1^{(1)}, k_1^{(4)} > 3r\right), d\left(k_1^{(2)}, k_1^{(3)} > 3r\right)\Big) \right\}$$

$$\times \frac{1}{(1+B^{j_1}d(k_1^{(1)}, k_1^{(2)}))^M} \frac{1}{(1+B^{j_1}d(k_1^{(3)}, k_1^{(4)}))^M} \frac{1}{(1+B^{j_2}d(k_2^{(1)}, k_2^{(3)}))^M}$$

$$\times \frac{1}{(1+B^{j_2}d(k_2^{(2)}, k_2^{(4)}))^M} \frac{1}{(1+B^{j_3}d(k_3^{(1)}, k_3^{(4)}))^M} \frac{1}{(1+B^{j_3}d(k_3^{(2)}, k_3^{(3)}))^M}$$

$$\leq \frac{C}{B^{4j_3}} \prod_{i=1}^4 \sum_{k_3^i} \left\{ \mathbb{I}\Big(d\left(k_1^{(1)}, k_1^{(4)} > 3r\right)\Big) \right\} \frac{1}{(1+B^{j_3}d(k_3^{(1)}, k_3^{(4)}))^M}$$

$$\times \frac{1}{(1+B^{j_3}d(k_3^{(2)}, k_3^{(3)}))^M}$$

$$\leq \frac{CB^{2j_3}}{B^{4j_3}} \sum_{k_3^{(1)}, k_3^{(4)}} \mathbb{I}\left(d(k_1^{(1)}, k_1^{(4)} > 3r)\right) \frac{1}{(1 + B^{j_3}d(k_3^{(1)}, k_3^{(4)}))^M}$$

$$\leq \frac{C}{B^{2j_3}} B^{\frac{1}{2}j_1 - j_3} B^{2j_3} = O(B^{\frac{1}{2}j_1 - j_3}) = O(B^{-\frac{1}{2}j_1})$$

because

$$\sum_{k' \in \mathcal{Z}_{j_3}} \mathbb{I}(d(k, k') \;\; > \;\; r) \frac{1}{(1 + B^{j_3}d(k, k'))^M} \leq CB^{2j_3} \int_r^\pi \frac{\sin\theta}{(1 + B^{j_3}\theta)^M} d\theta$$

$$\leq \;\; CB^{2j_3} \int_r^\pi \frac{\theta d\theta}{B^{Mj_3}\theta^M} \leq CB^{(M-2)(\frac{1}{2}j_1 - j_3)},$$

by taking $r = B^{-j_1/2}$, compare *Lemma 10* in [6]. The argument for the second and third term is analogous. Concerning the last summand, we recall that

$$Card\{k' \in \mathcal{V}_{j_1/2}(x_i) \cap \mathcal{X}_{j_3}\} \approx B^{2(j_3 - j_1/2)}, \tag{21}$$

for every $k \in \mathcal{X}_{j_1/2}$. Now denote

$$\Omega(k; j_3) \;\; := \;\; \{(k_1, k_2, k_3, k_4) : k_1, k_2, k_3, k_4 \in \mathcal{X}_{j_3}, k_1 \in \mathcal{V}_{j_3}(k), d(k_1, k_4)$$
$$\leq 3r, d(k_2, k_3) \leq 3r\},$$

where $r = \kappa B^{-\frac{1}{2}j_1}$. Heuristically, the idea is to split $\Omega(k)$ into regions where $(k_1, k_2, k_3, k_4)$ are each "close" to all three others, and regions where they are close two by two but the two pairs need not belong to the same neighbourhood. More precisely, we take $\Omega(k; j_3) \subseteq \Delta_1 \cup \Delta_2$, where

$$\Delta_1 = \{(k_1, k_2, k_3, k_4) : k_1, \ldots, k_4 \in B_{6r}(k), k_1 \in \mathcal{V}_{j_1/2}(k)\}$$

and

$$\Delta_2 = \{(k_1, k_2, k_3, k_4) : k_1 \in \mathcal{V}_{j_1/2}(k), \;\; k_2 \in S^2/B_{6r}(k)\}.$$

We can hence define

$$\Omega_1(k; j_3) := \Omega(k; j_3) \cap \Delta_1, \;\; \Omega_2(k; j_3) := \Omega(k; j_3) \cap \Delta_2.$$

In the region $\Omega_1(k; j_3)$, the idea is to keep $k$ fixed and then proceed by evaluating the cardinality of $B_{6r}(k)$; in view of (21), this leads to

$$(20) \leq \frac{C}{B^{4j_3}} \sum_{k \in \mathcal{X}_{j_1/2}} \sum_{(k_1, k_2, k_3, k_4) \in \Omega_1(k; j_3),} \frac{1}{(1 + B^{j_3}d(k_1, k_4))^M} \frac{1}{(1 + B^{j_3}d(k_2, k_3))^M}$$

$$= O\left(\frac{B^{4(j_3 - j_1/2) + j_1}}{B^{4j_3}}\right) = O(B^{-j_1}),$$

by (9) and *Lemma 4.8* in [31]. On the other hand, in the region $\Omega_2(k)$, we exploit the fact that $d(k_1, k_2), d(k_3, k_4) \geq 3r = 3\kappa B^{-j_1/2}$, so that we obtain the upper bound (for some $C > 0$)

$$
(20) \leq \frac{C}{B^{4j_3}} \sum_{k \in \mathcal{X}_{j_1/2}} \sum_{(k_3^{(1)}, \ldots, k_3^{(4)}) \in \Omega_2(k;j_3)} \frac{1}{(1 + B^{j_2} d(k_2^{(1)}, k_2^{(2)}))^M}
$$

$$
\times \frac{1}{(1 + B^{j_2} d(k_2^{(3)}, k_2^{(4)}))^M} \frac{1}{(1 + B^{j_3} d(k_3^{(1)}, k_3^{(4)}))^M} \frac{1}{(1 + B^{j_3} d(k_3^{(2)}, k_3^{(3)}))^M}
$$

$$
\leq \frac{C}{B^{4j_3}} \sum_{k \in \mathcal{X}_{j_1}} \sum_{(k_1, k_2, k_3, k_4) \in \Omega_2(k;j_3)} \frac{1}{(1 + C(B^{j_2 - j_1/2})^{2M})} \frac{1}{(1 + B^{j_3} d(k_1, k_4))^M}
$$

$$
\times \frac{1}{(1 + B^{j_3} d(k_2, k_3))^M}
$$

$$
= O\left( \frac{B^{4j_3 - 2M(j_2 - j_1/2)}}{B^{4j_3}} \right) = O(B^{-Mj_1}).
$$

The proof for the remaining terms is very similar and hence omitted for brevity's sake. ☐

### *3.2. Unknown angular power spectrum*

As the final result of this Section, we wish to focus on the case where the variance of the needlets coefficients is unknown and estimated from the data. A natural estimator for $\sigma_j$ is provided by

$$
\widetilde{\sigma}_j^2 = \frac{1}{N_j} \sum_{\xi_{jk} \in \mathcal{X}_j} |\beta_{jk}|^2
$$

where as before $N_j = card\{\mathcal{X}_j\} \approx B^{2j}$. We define our studentized statistics as $\widetilde{\beta}_{jk} := \beta_{jk}/\widetilde{\sigma}_j$ and we then consider

$$
\widetilde{I}_{j_1 j_2 j_3} = \sum_{k_3} \widetilde{\beta}_{j_1 k_1} \widetilde{\beta}_{j_2 k_2} \widetilde{\beta}_{j_3 k_3} \delta_{j_1 j_2 j_3}(k_1, k_2, k_3) h_{j_1 j_2 j_3}(k_1, k_2, k_3).
$$

Our next result shows that this studentization procedure has no effect on asymptotic behaviour.

**Theorem 5.** *As $j_1 \to \infty$, we have $\{\sigma_{j_1}^2 \sigma_{j_2}^2 \sigma_{j_3}^2\}^{-1} \widetilde{\sigma}_{j_1}^2 \widetilde{\sigma}_{j_2}^2 \widetilde{\sigma}_{j_3}^2 \longrightarrow 1$ in probability, and hence*

$$
\widetilde{I}_{j_1 j_2 j_3} \to_d N(0, 1).
$$

*Proof.* By a standard application of Slutzki's theorem, the weak convergence result follows immediately from the consistency of the variance estimator. We provide a proof for the three cases separately.

For case 1), i.e. $j_1 < j_2 < j_3$ we have immediately $E[\{\sigma_{j_1}^2 \sigma_{j_2}^2 \sigma_{j_3}^2\}^{-1} \widetilde{\sigma}_{j_1}^2 \widetilde{\sigma}_{j_2}^2 \widetilde{\sigma}_{j_3}^2] = 1$. Moreover

$$Var[\widetilde{\sigma}_{j_1}^2 \widetilde{\sigma}_{j_2}^2 \widetilde{\sigma}_{j_3}^2] = \frac{1}{N_{j_1}^2 N_{j_2}^2 N_{j_3}^2} \prod_{i=1}^{3} \sum_{\xi_{j_i k_i}, \xi_{j_i k_i'} \in \mathcal{X}_{j_i}} E[|\beta_{j_i k_i}|^2 |\beta_{j_i k_i'}|^2] - \sigma_{j_1}^4 \sigma_{j_2}^4 \sigma_{j_3}^4,$$

where

$$
\begin{aligned}
\sum_{\xi_{j_1 k_1}, \xi_{j_1 k_1'} \in \mathcal{X}_{j_1}} E[|\beta_{j_1 k_1}|^2 |\beta_{j_1 k_1'}|^2] &= \left( \sum_{\xi_{j_1 k_1} \in \mathcal{X}_{j_1}} E|\beta_{j_1 k_1}|^2 \right)^2 \\
&\quad + 2 \sum_{\xi_{j_1 k_1}, \xi_{j_1 k_1'} \in \mathcal{X}_{j_1}} |E[\beta_{j_1 k_1}, \beta_{j_1 k_1'}]|^2 \\
&= N_{j_1}^2 \sigma_{j_1}^4 + O\left( N_{j_1} \sigma_{j_1}^4 \right).
\end{aligned}
$$

Hence

$$
\begin{aligned}
Var[\widetilde{\sigma}_{j_1}^2 \widetilde{\sigma}_{j_2}^2 \widetilde{\sigma}_{j_3}^2] &= \frac{1}{N_{j_1}^2} \frac{1}{N_{j_2}^2} \frac{1}{N_{j_3}^2} (N_{j_1}^2 + O(N_{j_1}))(N_{j_2}^2 + O(N_{j_2})) \\
&\quad \times (N_{j_3}^2 + O(N_{j_3})) \sigma_{j_1}^4 \sigma_{j_2}^4 \sigma_{j_3}^4 - \sigma_{j_1}^4 \sigma_{j_2}^4 \sigma_{j_3}^4 \\
&= O\left( \frac{\sigma_{j_1}^4 \sigma_{j_2}^4 \sigma_{j_3}^4}{N_{j_1} N_{j_2} N_{j_3}} \right) = O\left( \frac{1}{N_{j_1} N_{j_2} N_{j_3}} \right).
\end{aligned}
$$

For case 2), we focus on the case where $j_1 = j_2 < j_3$; the remaining part of the proof is nearly identical. First note that

$$
\begin{aligned}
E\left\{ \frac{\widetilde{\sigma}_{j_1}^4 \widetilde{\sigma}_{j_2}^2}{\sigma_{j_1}^4 \sigma_{j_2}^2} \right\} &= \frac{\sum_{\xi_{j_2 k_2} \in \mathcal{X}_{j_2}} E|\beta_{j_2 k_2}|^2}{N_{j_2} \sigma_{j_1}^4 \sigma_{j_2}^2} \\
&\quad \times \left[ \left( \frac{4\pi}{N_{j_1}} \sum_{\xi_{j_1 k_1} \in \mathcal{X}_{j_1}} E|\beta_{j_1 k_1}|^2 \right)^2 \right. \\
&\quad \left. + \frac{2(4\pi)^2}{N_{j_1}^2} \sum_{\xi_{j_1 k_1}, \xi_{j_1 k_1'} \in \mathcal{X}_{j_1}} |E[\beta_{j_1 k_1}, \beta_{j_1 k_1'}]|^2 \right] \\
&= \frac{1}{\sigma_{j_1}^4 \sigma_{j_2}^2} (\sigma_{j_1}^4 + O\left( \frac{\sigma_{j_1}^4}{N_{j_1}} \right)) \sigma_{j_2}^2 = 1 + O\left( \frac{1}{N_{j_1}} \right).
\end{aligned}
$$

Likewise

$$
\begin{aligned}
Var[\widetilde{\sigma}_{j_1}^4 \widetilde{\sigma}_{j_2}^2] &= E\left( \frac{1}{N_{j_1}} \sum_{\xi_{j_1 k_1} \in \mathcal{X}_{j_1}} |\beta_{j_1 k_1}|^2 \right)^4 E\left( \frac{1}{N_{j_2}} \sum_{\xi_{j_2 k_2} \in \mathcal{X}_{j_2}} |\beta_{j_2 k_2}|^2 \right)^2 \\
&\quad - \left( \sigma_{j_1}^4 + O\left( \frac{\sigma_{j_1}^4}{N_{j_1}} \right) \right)^2 \sigma_{j_2}^4
\end{aligned}
$$

$$
= \frac{1}{N_{j_1}^4 N_{j_2}^2} E\left[\prod_{i=1}^4 \sum_{\xi_{j_1 k_1^i} \in \mathcal{X}_{j_1}} |\beta_{j_1 k_1^1}|^2 |\beta_{j_1 k_1^2}|^2 |\beta_{j_1 k_1^3}|^2 |\beta_{j_1 k_1^4}|^2\right]
$$

$$
\times E\left(\sum_{\xi_{j_2 k_2} \in \mathcal{X}_{j_2}} |\beta_{j_2 k_2}|^2\right)^2 - \left(\sigma_{j_1}^4 + O\left(\frac{\sigma_{j_1}^4}{N_{j_1}}\right)\right)^2 \sigma_{j_2}^4
$$

$$
= \left[\left(\sum_{\xi_{j_1 k_1} \in \mathcal{X}_{j_1}} E|\beta_{j_1 k_1}|^2\right)^4 + 12\left(\sum_{\xi_{j_1 k_1} \in \mathcal{X}_{j_1}} E|\beta_{j_1 k_1}|^2\right)^2\right.
$$

$$
\times \sum_{\xi_{j_1 k_1}, \xi_{j_1 k_1'} \in \mathcal{X}_{j_1}} |E[\beta_{j_1 k_1}, \beta_{j_1 k_1'}]|^2 + 8\left(\sum_{\xi_{j_1 k_1} \in \mathcal{X}_{j_1}} E|\beta_{j_1 k_1}|^2\right)
$$

$$
\times \left(\prod_{i=1}^3 \sum_{\xi_{j_1 k_1^i} \in \mathcal{X}_{j_1}} E[\beta_{j_1 k_1^1}, \beta_{j_1 k_1^2}] E[\beta_{j_1 k_1^2}, \beta_{j_1 k_1^3}] E[\beta_{j_1 k_1^1}, \beta_{j_1 k_1^3}]\right)
$$

$$
\times \frac{\sigma_{j_2}^4 (1 + O\left(\frac{1}{N_{j_2}}\right))}{N_{j_1}^4} - \left(\sigma_{j_1}^4 + O\left(\frac{\sigma_{j_1}^4}{N_{j_1}}\right)\right)^2 \sigma_{j_2}^4.
$$

By (9) and [31] (Lemma 4.8), we have

$$
\prod_{i=1}^3 \sum_{\xi_{j_1 k_1^i} \in \mathcal{X}_{j_1}} E[\beta_{j_1 k_1^1}, \beta_{j_1 k_1^2}] E[\beta_{j_1 k_1^2}, \beta_{j_1 k_1^3}] E[\beta_{j_1 k_1^1}, \beta_{j_1 k_1^3}]
$$

$$
\leq \sigma_{j_1}^6 \prod_{i=1}^3 \sum_{\xi_{j_1 k_1^i} \in \mathcal{X}_{j_1}} \frac{1}{(1 + B^{j_1} d(k_1^1, k_1^2))^3} \frac{1}{(1 + B^{j_1} d(k_1^2, k_1^3))^3} \frac{1}{(1 + B^{j_1} d(k_1^1, k_1^3))^3}
$$

$$
\leq \sigma_{j_1}^6 \sum_{\xi_{j_1 k_1^1}, \xi_{j_1 k_1^3} \in \mathcal{X}_{j_1}} \frac{C_M}{(1 + B^{j_2} d(k_1^1, k_1^3))^6} \leq C_M N_{j_1} \sigma_{j_1}^6.
$$

Hence the variance is bounded by

$$
\left| [N_{j_1}^4 \sigma_{j_1}^8 + 8\sigma_{j_1}^8 O\left(N_{j_1}^2\right) + 12\sigma_{j_1}^8 O\left(N_{j_1}^3\right)] \times \frac{\sigma_{j_2}^4 (1 + O\left(\frac{1}{N_{j_2}}\right))}{N_{j_1}^4} \right.
$$

$$
\left. - \left(\sigma_{j_1}^4 + O\left(\frac{\sigma_{j_1}^4}{N_{j_1}}\right)\right)^2 \sigma_{j_2}^4 \right| = O\left(\sigma_{j_1}^8 \sigma_{j_2}^4 \left(\frac{1}{N_{j_1}} + \frac{1}{N_{j_2}}\right)\right) = O\left(\frac{1}{N_{j_1}}\right).
$$

The proof for case 3) is an easy consequence of results by [6, 7]. □

## 4. Convergence to multiparameter Gaussian processes

Our aim in this Section is to extend the previous results to functional convergence theorems. The motivation for such an extension can be easily explained.

Indeed, from the applications points of view, practitioners are typically interested not only at the possible existence of non-Gaussianity and/or other features, but also to their location in frequency space. If we focus for instance on cosmological applications, which are the main driving rationale behind our work, it is important to recall that the existence of possible non-Gaussianities takes a very different meaning according to the scales where they are located, so that a suitable statistical procedure should provide information not only on their existence, but also on their position in the frequency domain. As an example, a huge debate has arisen in the Cosmological literature on the possible existence of a non-Gaussian "Cold Spot" in CMB data, much of the related literature concerning the determination of the angular scale of such features, see for instance [12, 13]. Concerning this feature, it may even be of interest to test for Gaussianity only on a subspace of the sphere (this is indeed what happens in practice, because of missing observations). The modification of (13) under these circumstances is straightforward: we would simply restrict our sum to a subset of the cubature points. Our following discussion would be asymptotically unaltered.

In [11, 23], it was proposed to build alternative forms of partial sum process from the bispectrum $B_{l_1 l_2 l_3}$, and to use them as a probe of non-Gaussianity at various scales. All different proposals were univariate, in the following sense. Assume the resolution of the experiment is such that frequencies up to $l = 1, \ldots, L$ are observed; the partial sums were then run only on a subset of configurations with cardinality $L$, whereas the number of multipole combinations $(l_1, l_2, l_3)$ which would be available for statistical analysis is in the order of $L^3$. One of the reasons for this restriction had to do with computational complexity: the evaluation of even a single bispectrum statistic is extremely time consuming, so that the exploration of all possible configurations is likely to be unfeasible even on the greatest supercomputing facilities. On the contrary, needlets are extremely convenient from a computational point of view, and there is no obstacle in considering larger frequency configurations, provided of course that the triangle conditions are satisfied.

We shall hence focus on two possible partial sums processes, which correspond broadly to cases 1 and 2 of the previous section; more precisely

$$J_{1L}(r_1, r_2) = \frac{1}{L} \sum_{j_1=1}^{[Lr_1]} \sum_{m=0}^{[Lr_2]} \widehat{I}_{j_1, j_1+K+m, j_1+K+m}, \tag{22}$$

$$J_{2L}(r) = \frac{1}{\sqrt{L}\sqrt{\sum_{m_1=0}^{N_1}(N(m_1)+1)}} \sum_{j=1}^{[Lr]} \sum_{m_1=0}^{N_1} \sum_{m_2=0}^{N(m_1)} \widehat{I}_{j_1, j_1+K+m_1, j_1+2K+m_1+m_2}, \tag{23}$$

where

$$\widehat{I}_{j_1, j_2, j_3} := \frac{I_{j_1, j_2, j_3}}{\sqrt{EI_{j_1, j_2, j_3}^2}},$$

$K \geq 0$ is some integer satisfying the constraint determined in Lemma 1,

$$N_1 = \max \left\{ m : 1 + B^{K+m} \geq B^{2K+m} \right\}$$

$$N(x) := \max \left\{ x : \log_B \left( 1 + B^{K+x} \right) - 2K - m_1 \geq 0 \right\}.$$

**Theorem 6.** *a) As $L \to \infty$*

$$J_{1L}(r_{1,r_2}) \Rightarrow W(r_1, r_2), \ in \ D[0,1]^2, \tag{24}$$

*where $W(.)$ is two-dimensional Brownian sheet, i.e. the zero mean Gaussian process with covariance function $EW(r_1, r_2)W(s_1, s_2) = (r_1 \wedge s_1)(r_2 \wedge s_2)$.*

*b) As $L \to \infty$*

$$J_{2L}(r) \Rightarrow W(r), \ in \ D[0,1] \tag{25}$$

*where $W(.)$ is standard Brownian Motion.*

*Here, $\Rightarrow$ denotes weak convergence in the sense of [9], $D[0,1]^p$, $p \in \mathbb{N}$ is the usual multidimensional Skorohod space.*

*Proof.* We start from (24); as usual, we need to establish convergence of the finite-dimensional distributions and tightness. Obviously,

$$EJ_{1L}(r_1, r_2) = 0,$$

and because $E[\hat{I}_{j_1,j_2,j_2}\hat{I}_{j_1',j_2',j_2'}] = \delta_{j_1}^{j_1'}\delta_{j_2}^{j_2'}$,

$$\begin{aligned}
&EJ_{1L}(r_1, r_2)J_{1L}(s_1, s_2) \\
&= \frac{1}{L^2} \sum_{j_1=1}^{[Lr_1]} \sum_{m=0}^{[Lr_2]} \sum_{j_1'=1}^{[Ls_1]} \sum_{m'=0}^{[Ls_2]} E\hat{I}_{j_1,j_1+K+m,j_1+K+m}\hat{I}_{j_1,j_1+K+m',j_1+K+m'} \\
&= \frac{1}{L^2} \sum_{j_1=1}^{[Lr_1]\wedge[Ls_1]} \sum_{m=0}^{[Lr_2]\wedge[Ls_2]} E\hat{I}_{j_1,j_1+K+m,j_1+K+m}^2 \\
&= \sum_{j_1=1}^{[Lr_1]\wedge[Ls_1]} \frac{[Lr_2] \wedge [Ls_2]}{L^2} \to [r_2 \wedge s_2][r_1 \wedge s_1].
\end{aligned}$$

To establish Gaussianity, we can again rely on the results by [32] and proceed with the bounds for the fourth-order cumulants. As the computations are very much the same as in the previous Section, we omit the details for brevity's sake. To consider tightness, we use the classical criteria given for instance in [35]. Define first the two-dimensional increments

$$J_{1L}((s_1, r_1] \times (s_2, r_2]) := J_{1L}(r_1, r_2) - J_{1L}(r_1, s_2) - J_{1L}(s_1, r_2) + J_{1L}(s_1, s_2).$$

It is again a standard computation to show that, as $L \to \infty$,

$$
\begin{aligned}
EJ_{1L}^2((s_1, r_1) \times (s_2, r_2)) &= \frac{1}{L^2} \bigg\{ \sum_{j_1=[Ls_1]+1}^{[Lr_1]} \sum_{m=[Ls_2]+1}^{[Lr_2]} E\hat{I}_{j_1,j_1+K+m,j_1+K+m}^2 \bigg\} \\
&= \frac{1}{L^2} \{[Lr_1] - [Ls_1]\}\{[Lr_2] - [Ls_2]\} \\
&\leq 4\{r_1 - s_1\}\{r_2 - s_2\}.
\end{aligned}
$$

We can then establish tightness by showing that

1) $s_i \leq t_i \leq r_i, i = 1, 2$

$$
\begin{aligned}
&E[J_{1L}^2((s_1, t_1) \times (s_2, t_2))J_{1L}^2((t_1, r_1) \times (t_2, r_2))] \\
&= \frac{1}{L^4} E\bigg( \sum_{j_1,j_1'=[Ls_1]}^{[Lt_1]} \sum_{\substack{j_2=j_1+K+m, \\ m=[Ls_2]+1}}^{[Lt_2]} \sum_{\substack{j_2'=j_1'+K+m', \\ m'=[Ls_2]+1}}^{[Lt_2]} \hat{I}_{j_1 j_2 j_2} \hat{I}_{j_1' j_2' j_2'} \bigg) \\
&\quad \times \bigg( \sum_{j_1,j_1'=[Lt_1]}^{[Lr_1]} \sum_{\substack{j_2=j_1+K+m, \\ m=[Lt_2]+1}}^{[Lr_2]} \sum_{\substack{j_2'=j_1'+K+m', \\ m'=[Lt_2]+1}}^{[Lr_2]} \hat{I}_{j_1 j_2 j_2} \hat{I}_{j_1' j_2' j_2'} \bigg) \\
&= \frac{1}{L^4} \bigg( \sum_{j_1=[Ls_1]}^{[Lt_1]} \sum_{m=[Ls_2]}^{[Lt_2]} E\hat{I}_{j_1,j_1+K+m,j_1+K+m}^2 \bigg) \\
&\quad \times \bigg( \sum_{j_1'=[Lt_1]}^{[Lr_1]} \sum_{m'=[Lt_2]}^{[Lr_2]} E\hat{I}_{j_1',j_1'+K+m',j_1'+K+m'}^2 \bigg) \\
&\leq 16(r_1 - t_1)(r_2 - t_2)(t_1 - s_1)(t_2 - s_2) \\
&\leq 4(r_1 - s_1)^2(r_2 - s_2)^2.
\end{aligned}
$$

2) $s_1 \leq r_1, s_2 \leq t_2 \leq r_2$,

$$
\begin{aligned}
&E[J_{1L}((s_1, r_1) \times (s_2, t_2))J_{1L}((s_1, r_1) \times (t_2, r_2))]^2 \\
&= \frac{1}{L^4} E\bigg( \sum_{j_1,j_1'=[Ls_1]}^{[Lt_1]} \sum_{\substack{j_2=j_1+K+m, \\ m=[Ls_2]+1}}^{[Lt_2]} \sum_{\substack{j_2'=j_1'+K+m', \\ m'=[Ls_2]+1}}^{[Lt_2]} \hat{I}_{j_1 j_2 j_2} \hat{I}_{j_1' j_2' j_2'} \bigg) \\
&\quad \times \bigg( \sum_{j_1,j_1'=[Ls_1]}^{[Lr_1]} \sum_{\substack{j_2=j_1+K+m, \\ m=[Lt_2]+1}}^{[Lr_2]} \sum_{\substack{j_2'=j_1'+K+m', \\ m'=[Lt_2]+1}}^{[Lr_2]} \hat{I}_{j_1 j_2 j_2} \hat{I}_{j_1' j_2' j_2'} \bigg) \\
&= \frac{1}{L^4} \bigg( \sum_{j_1=[Ls_1]}^{[Lr_1]} \sum_{m=[Ls_2]}^{[Lt_2]} E\hat{I}_{j_1,j_1+K+m,j_1+K+m}^2 \bigg) \\
&\quad \times \bigg( \sum_{j_1'=[Ls_1]}^{[Lr_1]} \sum_{m'=[Lt_2]}^{[Lr_2]} E\hat{I}_{j_1',j_1'+K+m',j_1'+K+m'}^2 \bigg)
\end{aligned}
$$

$$+ \frac{1}{L^4} \sum_{\substack{j_1^{(i)}=[Ls_1]+1, \\ i=1,..,4}}^{[Lr_1]} \sum_{\substack{m^{(i)}=[Ls_2]; \\ i=1,2}}^{[Lt_2]} \sum_{\substack{m^{(i)}=[Lt_2]; \\ i=3,4}}^{[Lr_2]}$$

$$\times \sum_{\gamma \in \Gamma_C} \rho\left(\gamma; \prod_{i=1}^{4} \widehat{I}_{j_1^{(i)},j_1^{(i)}+K+m^{(i)},j_1^{(i)}+K+m^{(i)}}\right)$$

For the first part it is easy to see that it is bounded by $16(r_1 - s_1)^2(r_2 - t_2)(t_2 - s_2)$; for the second part, for $j_1^{(i)}, j_2^{(i)}$ in each of their domain, we have,

$$\sum_{\gamma \in \Gamma_C} \rho\left(\gamma; \prod_{i=1}^{4} \widehat{I}_{j_1^{(i)} j_2^{(i)} j_2^{(i)}}\right)$$

$$= \sum_{\gamma \in \Gamma_{1C}} \rho(\gamma; \widehat{I}_{j_1 j_2^{(1)} j_2^{(1)}}^2 \widehat{I}_{j_1 j_2^{(3)} j_2^{(3)}}^2) \prod_{i=1}^{4} \delta_{j_1}^{j_1^{(i)}} \delta_{j_2}^{j_2^{(2)}} \delta_{j_2}^{j_2^{(4)}}_{j_2^{(3)}}$$

$$+ \sum_{\gamma \in \Gamma_{2C}} \rho(\gamma; \widehat{I}_{j_1 j_1' j_1'}^2 \widehat{I}_{j_1' j_2^{(3)} j_2^{(3)}}^2) \prod_{i=1}^{2} \delta_{j_2}^{j_1^{(i)}} \delta_{j_2}^{j_2^{(2)}} \delta_{j_2}^{j_2^{(4)}}_{j_2^{(3)}},$$

where $\Gamma_{1C}$ denotes the graphs with cliques (all nodes connected with all others), where $\Gamma_{2C}$ refers to graphs with loops of order two; these two disjoint classes cover all possible connected graphs with four nodes. We have

$$\frac{1}{L^4} \sum_{\substack{j_1^{(i)}=[Ls_1]+1, \\ i=1,..,4}}^{[Lr_1]} \sum_{\substack{m^{(i)}=[Ls_2]; \\ i=1,2}}^{[Lt_2]} \sum_{\substack{m^{(i)}=[Lt_2]; \\ i=3,4}}^{[Lr_2]} \sum_{\gamma \in \{fig2\}} \rho\left(\gamma; \prod_{i=1}^{4} \widehat{I}_{j_1^{(i)},j_1^{(i)}+K+m^{(i)},j_1^{(i)}+K+m^{(i)}}\right)$$

$$\leq \frac{C}{L^4} \sum_{j_1,j_1'=[Ls_1]+1}^{[Lr_1]} \sum_{\substack{m^{(i)}=[Ls_2]; \\ i=1,2}}^{[Lt_2]} \sum_{\substack{m^{(i)}=[Lt_2]; \\ i=3,4}}^{[Lr_2]} 1 \leq C'(r_1 - s_1)^2(r_2 - t_2)(t_2 - s_2)$$

from which we obtain

$$E[J_{1L}((s_1, r_1) \times (s_2, t_2))J_{1L}((s_1, r_1) \times (t_2, r_2))]^2 \leq C(r_1 - s_1)^2(r_2 - t_2)(t_2 - s_2)$$

Similarly, we can get the same result for
3) $s_1 \leq t_1 \leq r_1, s_2 \leq r_2$, that is

$$E[J_{1L}((s_1, t_1) \times (s_2, r_2))J_{1L}((t_1, r_1) \times (s_2, r_2))]^2 \leq C(r_2 - s_2)^2(r_1 - t_1)(t_1 - s_1)$$

This concludes the proof of (24).

For (25), we start again from the convergence of the finite-dimensional distributions; for notational simplicity, we stick to the univariate case. It is obvious that $EJ_{2L}(r) = 0$; on the other hand,

$$EJ_{2L}(r)J_{2L}(s)$$

$$= \sum_{j_1=1}^{[Lr]} \sum_{j_1'=1}^{[Ls]} \left\{ \sum_{m_1,m_1'=0}^{N_1} \sum_{m_2=0}^{N(m_1)} \sum_{m_2'=0}^{N(m_1')} \right\}$$

$$\times \frac{E\widehat{I}_{j_1,j_1+K+m_1,j_1+2K+m_1+m_2}\widehat{I}_{j_1',j_1'+K+m_1',j_1'+2K+m_1'+m_2'}}{L\sum_{m_1=0}^{N_1}(N(m_1)+1)}$$

$$= \frac{1}{L\sum_{m_1=0}^{N_1}(N(m_1)+1)} \sum_{j_1=1}^{[Lr]\wedge[Ls]} \sum_{m_1=0}^{N_1} \sum_{m_2=0}^{N(m_1)} E\widehat{I}_{j_1,j_1+K+m_1,j_1+2K+m_1+m_2}^2$$

$$= \frac{1}{L}\left([Lr]\wedge[Ls]-1\right) \to [r\wedge s], \text{ if } \min\{r,s\} > 0.$$

For Gaussianity, we analyze once again fourth-order cumulants, i.e. the connected components in the expansion of the fourth moment. As before, we need only focus on connected diagrams with four nodes, which can be partitioned into two classes: the cliques, where all nodes are connected with all three others, and diagrams with a loop of order 2. As before, these terms can be bounded by

$$\rho\left(\gamma; \prod_{i=1}^{4} \widehat{I}_{j_1^{(i)}j_2^{(i)}j_3^{(i)}}\right) = O(B^{-\max\{j_1^{(i)}\}/2}),$$

because for instance

$$\frac{1}{B^{j_3^{(1)}+j_3^{(2)}+2j_3^{(3)}}} \sum_{k_3^{(1)},k_3^{(2)},k_3^{(3)},k_3^{(4)}} \frac{1}{(1+B^{j_1^{(1)}}d(k_1^{(1)},k_1^{(2)}))^M} \frac{1}{(1+B^{j_2^{(1)}}d(k_2^{(1)},k_1^{(4)}))^M}$$

$$\times \frac{1}{(1+B^{j_2^{(2)}}d(k_2^{(2)},k_1^{(3)}))^M} \frac{1}{(1+B^{j_3^{(1)}}d(k_3^{(1)},k_2^{(3)}))^M} \frac{1}{(1+B^{j_3^{(2)}}d(k_3^{(2)},k_2^{(4)}))^M}$$

$$\times \frac{1}{(1+B^{j_3^{(3)}}d(k_3^{(3)},k_3^{(4)}))^M}$$

$$\leq \frac{\sum_{k_3^{(3)},k_3^{(4)}} \frac{1}{(1+B^{j_3^{(3)}}d(k_3^{(3)},k_3^{(4)}))^M}}{B^{j_3^{(1)}+j_3^{(2)}+2j_3^{3}}} \sum_{k_3^{(1)}} \frac{1}{(1+B^{j_3^{(1)}}d(k_3^{(1)},k_2^{(3)}))^M}$$

$$\times \sum_{k_3^{(2)}} \frac{1}{(1+B^{j_3^{(2)}}d(k_3^{(2)},k_2^{(4)}))^M}$$

$$\leq \frac{1}{B^{j_3^{(1)}+j_3^{(2)}+2j_3^{3}}} \sum_{k_3^{(3)},k_3^{(4)}} \frac{C_M}{(1+B^{j_3^{(3)}}d(k_3^{(3)},k_3^{(4)}))^M} \approx \frac{C_M B^{2j_3^{3}}}{B^{j_3^{(1)}+j_3^{(2)}+2j_3^{3}}}$$

$$= O(B^{-(j_3^{(1)}+j_3^{(2)})}).$$

To sum up

$$\prod_{i=1}^{4}\left\{\sum_{j_1^{(i)}=1}^{[Lr]}\sum_{m_1^{(i)}=0}^{N_1}\sum_{m_2^{(i)}=0}^{N\left(m_1^{(i)}\right)}\right\}\sum_{\gamma\in\Gamma_C}\rho\left(\gamma;\prod_{i=1}^{4}\widehat{I}_{j_1^{(i)},j_1^{(i)}+K+m_1^{(i)},j_1^{(i)}+2K+m_1^{(i)}+m_2^{(i)}}\right)$$

$$\leq\ C[Lr]\left(\sum_{m_1=0}^{N_1}\left(N\left(m_1\right)+1\right)\right)^2\sum_{j_1=1}^{[Lr]}B^{-j_1/2}=O(L).$$

It follows easily that

$$
\begin{aligned}
EJ_L^4(r) &= \sum_{j_1=1}^{[Lr]}\sum_{m_1=0}^{N_1}\sum_{m_2=0}^{N(m_1)}\frac{3\left(E\widehat{I}_{j_1,j_1+K+m_1,j_1+2K+m_1+m_2}^2\right)^2}{L^2\left(\sum_{m_1=0}^{N_1}(N(m_1)+1)\right)^2}\\
&\quad+\prod_{i=1}^{4}\left\{\sum_{j_1^{(i)}=1}^{[Lr]}\sum_{m_1^{(i)}=0}^{N_1}\sum_{m_2^{(i)}=0}^{N\left(m_1^{(i)}\right)}\right\}\\
&\quad\times\frac{\sum_{\gamma\in\Gamma_C}\rho(\gamma;\prod_{i=1}^{4}\widehat{I}_{j_1^{(i)},j_1^{(i)}+K+m_1^{(i)},j_1^{(i)}+2K+m_1^{(i)}+m_2^{(i)}})}{L^2\left(\sum_{m_1=0}^{N_1}(N(m_1)+1)\right)^2}\\
&= 3|EJ_L^2(r)|^2+O\left(\frac{1}{L}\right),
\end{aligned}
$$

which is enough to conclude the proof for the finite-dimensional distributions, in view of the standard argument from [32] that we used before.

To conclude the proof, we need only to consider tightness in $D([0,1])$. Note that $E[\widehat{I}_{j_1,j_2,j_3}\widehat{I}_{j_1',j_2',j_3'}]=\delta_{j_1}^{j_1'}\delta_{j_2}^{j_2'}\delta_{j_3}^{j_3'}$, the variance of $J_{2L}$ in $[s,r]$ is provided by

$$\mathbb{E}|J_{2L}(r)-J_{2L}(s)|^2=\sum_{j_1=[Ls]+1}^{[Lr]}\sum_{m_1=0}^{N_1}\sum_{m_2=0}^{N(m_1)}\frac{\mathbb{E}\widehat{I}_{j_1,j_1+K+m_1,j_1+2K+m_1+m_2}^2}{L\sum_{m_1=0}^{N_1}(N(m_1)+1)}\leq2(r-s)$$

Now we establish our tightness criterion. For any $0\leq s\leq t\leq r\leq1$,

$$
\begin{aligned}
&\mathbb{E}|J_{2L}(r)-J_{2L}(t)|^2|J_{2L}(t)-J_{2L}(s)|^2\\
&=\sum_{j_1=[Lt]}^{[Lr]}\sum_{m_1=0}^{N_1}\sum_{m_2=0}^{N(m_1)}\frac{\mathbb{E}\widehat{I}_{j_1,j_1+K+m_1,j_1+2K+m_1+m_2}^2}{L\sum_{m_1=0}^{N_1}(N(m_1)+1)}\\
&\quad\times\sum_{j_1'=[Ls]}^{[Lt]}\sum_{m_1'=0}^{N_1}\sum_{m_2'=0}^{N(m_1')}\frac{\mathbb{E}\widehat{I}_{j_1',j_1'+K+m_1',j_1'+2K+m_1'+m_2'}^2}{L\sum_{m_1=0}^{N_1}(N(m_1)+1)}\\
&\quad+\sum_{\substack{j_1^{(i)}=[Lt],\,j_1^{(i)}=[Ls],\\i=1,2\quad\quad i=3,4}}^{[Lr]}\sum^{[Lt]}\left\{\prod_{l=1}^{4}\sum_{m_1^{(l)}=0}^{N_1}\sum_{m_2^{(l)}=0}^{N(m_1^{(l)})}\right\}
\end{aligned}
$$

$$\times \sum_{\gamma \in \Gamma_C} \frac{\rho(\gamma; \prod_{l=1}^{4} \widehat{I}_{j_1^{(l)}, j_1^{(l)} + K + m_1^{(l)}, j_1^{(l)} + 2K + m_1^{(l)} + m_2^{(l)}})}{\left(L \sum_{m_1=0}^{N_1} (N(m_1) + 1)\right)^2}$$

$$\leq \frac{2}{L^2}([Lr] - [Lt])([Lt] - [Ls]) \leq 4(r - s)^2.$$

Thus we finished the proof of tightness. $\qquad \square$

## 5. Behaviour under non-Gaussianity

In this final Section, we shall provide some quick and informal discussion on the behaviour of our statistics under non-Gaussianity; see [33, 26] for other applications of the needlets to cosmological data analysis. There exist of course a huge variety of non-Gaussian models for spherical random fields, and we shall delay a much more detailed treatment to future work. Our purpose here is different, i.e. we want to provide some heuristic discussion on the expected behaviour of our procedures for physically motivated non-Gaussian models. This will provide some guidance to practitioners for applications to CMB data, which are currently under way in a separate work.

We start from the expected value of the needlets bispectrum, which is provided by

$$\mathbb{E}I_{j_1 j_2 j_3} \approx \frac{1}{\sqrt{B^{2j_3}}} \sum_{k_1 k_2 k_3} \mathbb{E}\widehat{\beta}_{j_1 k_1} \widehat{\beta}_{j_2 k_2} \widehat{\beta}_{j_3 k_3} \delta_{j_1 j_2 j_3}(k_1, k_2, k_3)$$

$$\approx \frac{1}{\sigma_{j_1}\sigma_{j_2}\sigma_{j_3}} \frac{1}{\sqrt{B^{2j_3}}} \sum_{k_1 k_2 k_3} \sum_{\substack{l_1, l_2, l_3 \\ m_1 m_2 m_3}} b\left(\frac{l_1}{B^{j_1}}\right) b\left(\frac{l_2}{B^{j_2}}\right)$$

$$\times b\left(\frac{l_3}{B^{j_3}}\right) \mathbb{E}\left(a_{l_1 m_1} a_{l_2 m_2} a_{l_3 m_3}\right)$$

$$\times Y_{l_1 m_1}(\xi_{j_1 k_1}) Y_{l_2 m_2}(\xi_{j_2 k_2}) Y_{l_3 m_3}(\xi_{j_3 k_3}) \delta_{j_1 j_2 j_3}(k_1, k_2, k_3)$$

$$= \frac{1}{\sigma_{j_1}\sigma_{j_2}\sigma_{j_3}} B^{j_3} \sum_{l_1, l_2, l_3 = B^{j-1}}^{B^{j+1}} \sum_{m_1 m_2 m_3} b\left(\frac{l_1}{B^{j_1}}\right) b\left(\frac{l_2}{B^{j_2}}\right) b\left(\frac{l_3}{B^{j_3}}\right) b_{l_1 l_2 l_3}$$

$$\times \begin{pmatrix} l_1 & l_2 & l_3 \\ m_1 & m_2 & m_3 \end{pmatrix} \begin{pmatrix} l_1 & l_2 & l_3 \\ 0 & 0 & 0 \end{pmatrix} \sqrt{\frac{(2l_1 + 1)(2l_2 + 1)(2l_3 + 1)}{4\pi}}$$

$$\times \left\{ \frac{1}{B^{2j_3}} \sum_{k_1 k_2 k_3} Y_{l_1 m_1}(\xi_{j_1 k_1}) Y_{l_2 m_2}(\xi_{j_2 k_2}) Y_{l_3 m_3}(\xi_{j_3 k_3}) \delta_{j_1 j_2 j_3}(k_1, k_2, k_3) \right\}.$$

Here, we recall that $b_{l_1 l_2 l_3}$ is the so-called reduced bispectrum (see for instance [20, 23]), which collects the non-Gaussian component in the third order moment

$\mathbb{E}\left(a_{l_1 m_1} a_{l_2 m_2} a_{l_3 m_3}\right)$, i.e. by definition

$$\mathbb{E}\left(a_{l_1 m_1} a_{l_2 m_2} a_{l_3 m_3}\right)$$
$$= \begin{pmatrix} l_1 & l_2 & l_3 \\ m_1 & m_2 & m_3 \end{pmatrix} \begin{pmatrix} l_1 & l_2 & l_3 \\ 0 & 0 & 0 \end{pmatrix} \sqrt{\frac{(2l_1+1)(2l_2+1)(2l_3+1)}{4\pi}} b_{l_1 l_2 l_3}.$$

In the cosmological literature, a very popular model for $b_{l_1 l_2 l_3}$ is provided by the so-called Sachs-Wolfe bispectrum ([20], equation (21)), which yields

$$b_{l_1 l_2 l_3} = -6 f_{NL} \left\{ C_{l_1} C_{l_2} + C_{l_1} C_{l_3} + C_{l_2} C_{l_3} \right\},$$

where $f_{NL}$ is a physical constant (see for instance [8, 37]). The Wigner coefficients on the right hand side ensure that the expected value $\mathbb{E}\left(a_{l_1 m_1} a_{l_2 m_2} a_{l_3 m_3}\right)$ is rotationally invariant under a change of coordinate, an obvious consequence of the isotropy of the random field. For our purposes below, it is sufficient to recall that

$$\begin{pmatrix} l & l & l \\ 0 & 0 & 0 \end{pmatrix} \approx \frac{(-1)^{-3l/2}}{l}, \quad \begin{pmatrix} l_0 & l & l+l_0 \\ 0 & 0 & 0 \end{pmatrix} \approx \frac{(-1)^{-l_0+l}}{\sqrt{l}}. \quad (26)$$

Now exploiting again the cubature formula (3) we obtain as before

$$\left\{ \frac{1}{B^{2j_3}} \sum_{k_1 k_2 k_3} Y_{l_1 m_1}(\xi_{j_1 k_1}) Y_{l_2 m_2}(\xi_{j_2 k_2}) Y_{l_3 m_3}(\xi_{j_3 k_3}) \delta_{j_1 j_2 j_3}(k_1, k_2, k_3) \right\}$$
$$\approx \int_{\mathbb{S}^2} Y_{l_1 m_1}(\xi) Y_{l_2 m_2}(\xi) Y_{l_3 m_3}(\xi) d\xi$$
$$= \begin{pmatrix} l_1 & l_2 & l_3 \\ m_1 & m_2 & m_3 \end{pmatrix} \begin{pmatrix} l_1 & l_2 & l_3 \\ 0 & 0 & 0 \end{pmatrix} \sqrt{\frac{(2l_1+1)(2l_2+1)(2l_3+1)}{4\pi}},$$

whence

$$\mathbb{E} I_{j_1 j_2 j_3} = \frac{1}{\sigma_{j_1} \sigma_{j_2} \sigma_{j_3}} B^{j_3} \sum_{l_1, l_2, l_3 = B^{j-1}}^{B^{j+1}} b\left(\frac{l_1}{B^{j_1}}\right) b\left(\frac{l_2}{B^{j_2}}\right) b\left(\frac{l_3}{B^{j_3}}\right) b_{l_1 l_2 l_3} \begin{pmatrix} l_1 & l_2 & l_3 \\ 0 & 0 & 0 \end{pmatrix}^2$$
$$\times \frac{(2l_1+1)(2l_2+1)(2l_3+1)}{4\pi} \sum_{m_1 m_2 m_3} \begin{pmatrix} l_1 & l_2 & l_3 \\ m_1 & m_2 & m_3 \end{pmatrix}^2$$
$$= \frac{B^{j_3}}{\sigma_{j_1} \sigma_{j_2} \sigma_{j_3}} \sum_{l_1, l_2, l_3 = B^{j-1}}^{B^{j+1}} b\left(\frac{l_1}{B^{j_1}}\right) b\left(\frac{l_2}{B^{j_2}}\right) b\left(\frac{l_3}{B^{j_3}}\right) b_{l_1 l_2 l_3} \begin{pmatrix} l_1 & l_2 & l_3 \\ 0 & 0 & 0 \end{pmatrix}^2$$
$$\times \frac{(2l_1+1)(2l_2+1)(2l_3+1)}{4\pi},$$

in view of the orthonormality properties of the Wigner's $3j$ coefficients. To keep the analogy with the cosmological literature, we shall focus on "equilateral" and

"squeezed" configurations, see [3, 23]. In the equilateral case $j_1 = j_2 = j_3 = j$ we have

$$
\begin{aligned}
\mathbb{E} I_{jjj} \;=\; & \frac{B^j}{\sigma_{j_1}\sigma_{j_2}\sigma_{j_3}} \sum_{l_1,l_2,l_3=B^{j-1}}^{B^{j+1}} b\left(\frac{l_1}{B^j}\right) b\left(\frac{l_2}{B^j}\right) b\left(\frac{l_3}{B^j}\right) b_{l_1 l_2 l_3} \\
& \times \left(\begin{array}{ccc} l_1 & l_2 & l_3 \\ 0 & 0 & 0 \end{array}\right)^2 \frac{(2l_1+1)(2l_2+1)(2l_3+1)}{4\pi}.
\end{aligned}
$$

Now recall $l \approx B^j$, $\sigma_j^3 \approx C_{B^j}^{3/2} B^{3j}$, $b_{l_1 l_2 l_3} \approx f_{NL} l_1^{-\alpha} l_2^{-\alpha}$ so that, using also (17)

$$
\begin{aligned}
B^j & \sum_{l_1,l_2,l_3=B^{j-1}}^{B^{j+1}} b\left(\frac{l_1}{B^j}\right) b\left(\frac{l_2}{B^j}\right) b\left(\frac{l_3}{B^j}\right) \frac{b_{l_1 l_2 l_3}}{\sigma_j^3} \left(\begin{array}{ccc} l_1 & l_2 & l_3 \\ 0 & 0 & 0 \end{array}\right)^2 \\
& \times \frac{(2l_1+1)(2l_2+1)(2l_3+1)}{4\pi} \\
\simeq \; & B^j \sum_{l_1,l_2,l_3=B^{j-1}}^{B^{j+1}} b\left(\frac{l_1}{B^j}\right) b\left(\frac{l_2}{B^j}\right) b\left(\frac{l_3}{B^j}\right) \frac{f_{NL} l_1^{-\alpha} l_2^{-\alpha}}{C_{B^j}^{3/2} B^{3j}} \left(\begin{array}{ccc} l_1 & l_2 & l_3 \\ 0 & 0 & 0 \end{array}\right)^2 B^{3j} \\
\simeq \; & B^j \sum_{l_1,l_2,l_3=B^{j-1}}^{B^{j+1}} b\left(\frac{l_1}{B^j}\right) b\left(\frac{l_2}{B^j}\right) b\left(\frac{l_3}{B^j}\right) \frac{f_{NL} B^{-2j\alpha}}{B^{-3j\alpha/2} B^{3j}} B^j \simeq f_{NL} B^{2j} B^{-j\alpha/2}.
\end{aligned}
$$

This suggests the expected value of the needlets bispectrum can either diverge or converge to zero, according to the asymptotic behaviour of the angular power spectrum; in particular, it does diverge for all $\alpha < 4$. On the other hand, for $j_1 << j_2 = j_3$ by an analogous argument we obtain

$$
\begin{aligned}
\mathbb{E} I_{j_1 j_2 j_2} \simeq \; & \frac{f_{NL} B^{j_2}}{B^{-j_1(\alpha/2-1)} B^{-j_2(\alpha-2)}} \\
& \times \sum_{l_1,l_2,l_3=B^{j-1}}^{B^{j+1}} b\left(\frac{l_1}{B^{j_1}}\right) b\left(\frac{l_2}{B^{j_2}}\right) b\left(\frac{l_3}{B^{j_2}}\right) l_1^{-\alpha} l_2^{-\alpha} \left(\begin{array}{ccc} l_1 & l_2 & l_3 \\ 0 & 0 & 0 \end{array}\right)^2 B^{j_1} B^{2j_2} \\
\simeq \; & \frac{f_{NL} B^{j_2}}{B^{-j_1(\alpha/2-1)} B^{-j_2(\alpha-2)}} \\
& \times \sum_{l_1,l_2,l_3=B^{j-1}}^{B^{j+1}} b\left(\frac{l_1}{B^{j_1}}\right) b\left(\frac{l_2}{B^{j_2}}\right) b\left(\frac{l_3}{B^{j_2}}\right) B^{-j_1\alpha} B^{-j_2\alpha} B^{j_2} \\
\simeq \; & f_{NL} B^{-j_1\alpha/2} B^{2j_2}.
\end{aligned}
$$

As for the usual bispectrum, the previous computations suggest that the power is maximized by "squeezing" frequencies, i.e. maximizing the differences between the "side lengths" $j_1$ and $j_2$. This is the same sort of qualitative result which was found for the bispectrum in [23] and successfully applied to CMB data in [11]. Our heuristic calculations in this Section suggest very clearly that the needlets bispectrum may enjoy the same good power properties, at the same time healing the difficulties that were met in [11] in the presence of missing observations. These claims will be soon scrutinized on simulations and on real data.